\documentclass[11pt]{jdg-p}
\usepackage{amsfonts}
\usepackage{amssymb}
\usepackage{eucal}

\begin{document}

\font\eightrm=cmr8
\font\eightit=cmti8
\font\eighttt=cmtt8
\def\tci
{\hbox{\hskip1.8pt$\rightarrow$\hskip-11.5pt$^{^{C^\infty}}$\hskip-1.3pt}}
\def\nft
{\hbox{$n$\hskip3pt$\equiv$\hskip4pt$5$\hskip4.4pt$($mod\hskip2pt$3)$}}
\def\bbR{\text{\bf R}}
\def\rto{\mathbf{R}\hskip-.5pt^2}
\def\rtr{\text{\bf R}\hskip-.7pt^3}
\def\rfo{\text{\bf R}\hskip-.7pt^4}
\def\rn{{\mathbf{R}}^{\hskip-.6ptn}}
\def\bbRP{\text{\bf R}\text{\rm P}}
\def\bbC{\text{\bf C}}
\def\bbZ{\text{\bf Z}}
\def\hyp{\hskip.5pt\vbox
{\hbox{\vrule width3ptheight0.5ptdepth0pt}\vskip2.2pt}\hskip.5pt}
\def\er{r}
\def\es{s}
\def\df{d\hskip-.8ptf}
\def\fv{\mathcal{F}}
\def\fvp{\fv_{\nh p}}
\def\wv{\mathcal{W}}
\def\vt{\mathcal{P}}
\def\tv{\mathcal{T}}
\def\vtx{\vt_{\nh x}}
\def\fh{f}
\def\rc{\theta}
\def\jm{\mathcal{I}}
\def\ke{\mathcal{K}}
\def\xc{\mathcal{X}_c}
\def\lz{\mathcal{L}}
\def\lo{\lz_0}
\def\xe{\mathcal{E}}
\def\eo{\xe_0}
\def\hm{\hskip1.9pt\widehat{\hskip-1.9ptM\hskip-.2pt}\hskip.2pt}
\def\hmt{\hskip1.9pt\widehat{\hskip-1.9ptM\hskip-.5pt}_t}
\def\hmz{\hskip1.9pt\widehat{\hskip-1.9ptM\hskip-.5pt}_0}
\def\hmp{\hskip1.9pt\widehat{\hskip-1.9ptM\hskip-.5pt}_p}
\def\hg{\hskip1.2pt\widehat{\hskip-1.2ptg\hskip-.4pt}\hskip.4pt}
\def\hq{\hskip1.5pt\widehat{\hskip-1.5ptQ\hskip-.5pt}\hskip.5pt}
\def\q{q}
\def\bq{\hat q}
\def\p{p}
\def\w{\vt^\perp}
\def\x{v}
\def\y{y}
\def\vp{\vt^\perp}
\def\vd{\vt\hh'}
\def\vdx{\vd{}\hskip-4.5pt_x}
\def\bz{b\hh}
\def\fe{F}
\def\fy{\phi}
\def\vl{\Lambda}
\def\hy{\mathcal{V}}
\def\vh{h}
\def\mv{V}
\def\vo{V_{\nnh0}}
\def\ao{A_0}
\def\bo{B_0}
\def\uv{\mathcal{U}}
\def\sv{\mathcal{S}}
\def\svp{\sv_p}
\def\xv{\mathcal{X}}
\def\xvp{\xv_p}
\def\yv{\mathcal{Y}}
\def\yvp{\yv_p}
\def\zv{\mathcal{Z}}
\def\zvp{\zv_p}
\def\cv{\mathcal{C}}
\def\dy{\mathcal{D}}
\def\nv{\mathcal{N}}
\def\iv{\mathcal{I}}
\def\gkp{\Sigma}
\def\hs{\hskip.7pt}
\def\hh{\hskip.4pt}
\def\nh{\hskip-.7pt}
\def\nnh{\hskip-1pt}
\def\hrz{^{\hskip.5pt\text{\rm hrz}}}
\def\vrt{^{\hskip.2pt\text{\rm vrt}}}
\def\vt{\varTheta}
\def\vg{\varGamma}
\def\my{\mu}
\def\ny{\nu}
\def\gy{\lambda}
\def\ax{\alpha}
\def\bx{\beta}
\def\cx{\gamma}
\def\ay{a}
\def\by{b}
\def\cy{c}
\def\gp{\mathrm{G}}
\def\hp{\mathrm{H}}
\def\kp{\mathrm{K}}
\def\gm{\gamma}
\def\Gm{\Gamma}
\def\Lm{\Lambda}
\def\Dt{\Delta}
\def\sj{\sigma}
\def\lg{\langle}
\def\rg{\rangle}
\def\lr{\lg\,,\rg}
\def\vs{vector space}
\def\rvs{real vector space}
\def\vf{vector field}
\def\tf{tensor field}
\def\tvn{the vertical distribution}
\def\dn{distribution}
\def\pt{point}
\def\tc{tor\-sion\-free connection}
\def\ea{equi\-af\-fine}
\def\rt{Ric\-ci tensor}
\def\pde{partial differential equation}
\def\pf{projectively flat}
\def\pfs{projectively flat surface}
\def\pfc{projectively flat connection}
\def\pftc{projectively flat tor\-sion\-free connection}
\def\su{surface}
\def\sco{simply connected}
\def\psr{pseu\-\hbox{do\hs-}Riem\-ann\-i\-an}
\def\inv{-in\-var\-i\-ant}
\def\trinv{trans\-la\-tion\inv}
\def\feo{dif\-feo\-mor\-phism}
\def\feic{dif\-feo\-mor\-phic}
\def\feicly{dif\-feo\-mor\-phi\-cal\-ly}
\def\Feicly{Dif\-feo\-mor\-phi\-cal\-ly}
\def\diml{-di\-men\-sion\-al}
\def\prl{-par\-al\-lel}
\def\skc{skew-sym\-met\-ric}
\def\sky{skew-sym\-me\-try}
\def\Sky{Skew-sym\-me\-try}
\def\dbly{-dif\-fer\-en\-ti\-a\-bly}
\def\cs{con\-for\-mal\-ly symmetric}
\def\cf{con\-for\-mal\-ly flat}
\def\ls{locally symmetric}
\def\ecs{essentially con\-for\-mal\-ly symmetric}
\def\rr{Ric\-ci-re\-cur\-rent}
\def\kf{Killing field}
\def\om{\omega}
\def\vol{\varOmega}
\def\dv{\delta}
\def\ve{\varepsilon}
\def\zt{\zeta}
\def\kx{\kappa}
\def\mf{manifold}
\def\mfd{-man\-i\-fold}
\def\bmf{base manifold}
\def\bd{bundle}
\def\tbd{tangent bundle}
\def\ctb{cotangent bundle}
\def\bp{bundle projection}
\def\prc{pseu\-\hbox{do\hs-}Riem\-ann\-i\-an metric}
\def\prd{pseu\-\hbox{do\hs-}Riem\-ann\-i\-an manifold}
\def\Prd{pseu\-\hbox{do\hs-}Riem\-ann\-i\-an manifold}
\def\npd{null parallel distribution}
\def\pj{-pro\-ject\-a\-ble}
\def\pd{-pro\-ject\-ed}
\def\lcc{Le\-vi-Ci\-vi\-ta connection}
\def\vb{vector bundle}
\def\vbm{vec\-tor-bun\-dle morphism}
\def\kerd{\text{\rm Ker}\hskip2.7ptd}
\def\ro{\rho}
\def\sy{\sigma}
\def\ts{total space}
\def\pmb{\pi}

\newtheorem{theorem}{Theorem}[section] 
\newtheorem{proposition}[theorem]{Proposition} 
\newtheorem{lemma}[theorem]{Lemma} 
\newtheorem{corollary}[theorem]{Corollary} 
  
\theoremstyle{definition} 
  
\newtheorem{defn}[theorem]{Definition} 
\newtheorem{notation}[theorem]{Notation} 
\newtheorem{example}[theorem]{Example} 
\newtheorem{conj}[theorem]{Conjecture} 
\newtheorem{prob}[theorem]{Problem} 
  
\theoremstyle{remark} 
  
\newtheorem{remark}[theorem]{Remark}

\renewcommand{\thepart}{\Roman{part}}
\title[Compact manifolds with parallel Weyl tensor]{Compact 
pseu\-do\hs-Riem\-ann\-i\-an manifolds with parallel Weyl tensor}
\author[A. Derdzinski]{Andrzej Derdzinski} 
\address{Department of Mathematics, Ohio State University, 
Columbus, OH 43210} 
\email{andrzej@math.ohio-state.edu} 
\author[W. Roter]{Witold Roter}
\address{Institute of Mathematics and Computer Science, 
Wroc\l aw University of Technology, 
Wy\-brze\-\.ze Wys\-pia\'n\-skiego 27, 
50-370 Wroc\l aw, Poland}
\email{roter@im.pwr.wroc.pl} 
\def\leftmark{A. Derdzinski \&\ W. Roter}
\def\rightmark{Compact manifolds with parallel Weyl tensor}

\begin{abstract}
It is shown that in every dimension $\,n=3j+2$, $\hs j=1,2,3,\nh\dots\hs$, 
there exist compact pseu\-do\hs-Riem\-ann\-i\-an manifolds with parallel 
Weyl tensor, which are Ric\-ci-re\-cur\-rent, but neither con\-for\-mal\-ly 
flat nor locally symmetric, and represent all indefinite metric signatures. 
The manifolds in question are dif\-feo\-mor\-phic to nontrivial torus bundles 
over the circle. They all arise from a construction that a priori yields 
bundles over the circle, having as the fibre either a torus, or a $\,2$-step 
nil\-man\-i\-fold with a complete flat tor\-sion\-free connection; our 
argument only realizes the torus case.
\end{abstract}

\maketitle

\voffset=45pt\hoffset=57pt 

\section*{Introduction}
\setcounter{equation}{0}
A \prd\ $\,(M,g)\,$ of dimension $\,n\ge4\,$ is called {\em con\-for\-mal\-ly 
symmetric\/} \cite{chaki-gupta} if its Weyl conformal tensor is parallel. If, 
in addition, $\,(M,g)\,$ is neither con\-for\-mal\-ly flat nor locally 
symmetric, it is said to be {\em essentially con\-for\-mal\-ly symmetric}. 

All \ecs\ \prc s are indefinite \cite[Theorem~2]{derdzinski-roter-77}. 
Numerous examples of such metrics on open manifolds are known 
\cite{derdzinski-roter-07,roter}, which raises the question whether they exist 
on any compact manifolds, cf.\ \cite{simon}. This paper provides an answer:
\begin{theorem}\label{maith}In every dimension $\,n=3j+2\ge5$, 
$\,j=1,2,3,\dots\hs$, there exists a compact pseu\-do\hs-Riem\-ann\-i\-an 
manifold\/ $\,(M,g)\,$ of any prescribed indefinite metric signature, which is 
essentially con\-for\-mal\-ly symmetric, Ric\-ci-re\-cur\-rent, and 
dif\-feo\-mor\-phic to a torus bundle over the circle, but not homeomorphic 
to, or even covered by, the torus $\,T^{\hh n}\nnh$.
\end{theorem}
Here $\,(M,g)\,$ is called {\it Ric\-ci-re\-cur\-rent\/} if, for every tangent 
vector field $\,w$, the Ric\-ci tensor $\,\hs\text{\rm Ric}\hs\,$ and the 
covariant derivative $\,\nabla_{\!w}\hs\text{\rm Ric}\hs\,$ are linearly 
dependent at every point.

Each \mf\ in Theorem~\ref{maith} arises as the quotient $\,M=\hm\nnh/\hh\Gm\,$ 
for a suitable discrete group $\,\Gm\hs$ of isometries of its universal 
covering space $\,\hm$, \feic\ to $\,\rn\nnh$, with a metric belonging to a 
family constructed by the second author in \cite{roter}. For every dimension 
$\,n=3j+2$, the metrics in that family admitting such compact quotients form 
an in\-fi\-nite\diml\ space of local moduli (Remark~\ref{modul}). However, 
our argument provides no explicit descriptions of the metrics, or the groups 
$\,\Gm$.

Con\-for\-mal symmetry is one of the {\it natural linear conditions} in the 
sense of Besse \cite[p.\ 433]{besse} that can be imposed on the covariant 
derivatives of the irreducible components of the curvature tensor under the 
action of the pseu\-do\hs-or\-thog\-o\-nal group. The analogous conditions on 
the other two components characterize metrics having constant scalar curvature 
and, respectively, parallel Ric\-ci tensor, including the Einstein metrics. 
Compact Riemannian or K\"ahler manifolds of these two classes are the model 
cases of the Yamabe problem and Calabi's conjectures.

Compact con\-for\-mal\-ly symmetric manifolds have generated much less 
interest. However, those among them having the specific form 
$\,M=\hm\nnh/\hh\Gm\,$ mentioned above are related to another familiar class of 
geometric structures. Namely, we show, in Remarks~\ref{ntrst} and~\ref{nilbd}, 
that any such $\,M\,$ is a bundle over the circle, and its fibre is either a 
torus, or a $\,2$-step nil\-man\-i\-fold admitting a complete flat 
tor\-sion\-free connection with a nonzero parallel vector field. (Our argument 
only succeeds in realizing the torus case.) Complete flat tor\-sion\-free 
connections on compact manifolds are the subject of a vast literature, 
outlined in \cite{charette-drumm-goldman-morrill}, and on nil\-man\-i\-folds 
they exist relatively often, though not always \cite{benoist}.

One easily verifies that no \ecs\ \mf\ is locally reducible. The gaps in the 
dimension list of Theorem~\ref{maith} cannot therefore be filled with the aid 
of Riemannian products. Thus, Theorem~\ref{maith} leaves the existence 
question unanswered in dimensions $\,n\ge4\,$ other than those of the form 
$\,n=3j+2$. While for $\,n\ge5\,$ this may be due to the particular nature of 
our argument, designed to work only when $\,\nft$, the reason why 
Theorem~\ref{maith} fails to include the case $\,n=4\,$ seems less of a 
coincidence. In fact, using Theorem~\ref{nonex} of the present paper, we show 
in \cite{derdzinski-roter} that {\it every four\diml\ \ecs\ Lo\-rentz\-i\-an 
\mf\ is noncompact}.

There are further instances where particular details of Theorem~\ref{maith} 
reflect more general facts. Namely, two other results of 
\cite{derdzinski-roter} state that {\it the fundamental group of a compact 
\ecs\ \mf\ is always infinite}, and {\it for any compact \ecs\ 
Lo\-rentz\-i\-an \mf\/ $\,(M,g)$, some two-fold covering manifold of\/ 
$\,M\hs$ is a bundle over the circle and its fibre admits a flat 
tor\-sion\-free connection with a nonzero parallel vector field}. 

\section{Preliminaries}\label{prel}
Let a group $\,\Gm\hs$ act on a manifold $\,\hm\,$ freely by 
\feo s. The action of $\,\Gm\hs$ on $\,\hm\,$ is called {\it 
properly dis\-con\-tin\-u\-ous} if there exists a locally \feic\ surjective 
mapping $\,\pi:\hm\to M\,$ onto some manifold $\,M$ such that the 
$\,\pi$-pre\-im\-ages of points of $\,M\,$ are precisely the orbits of the 
$\,\Gm\hs$ action. (Cf.\ \cite[p.\ 187]{charette-drumm-goldman-morrill}.) We 
then refer to $\,M\,$ as the {\it quotient\/} of $\,\hm\,$ under the action of 
$\,\Gm\hs$ and write $\,M=\hm\nnh/\hh\Gm$.

The index $\,j=1,2,3,\dots\,$ is always used to label the terms of sequences, 
with $\,x_j\to x\,$ meaning that $\,x=\lim_{\,j\to\infty}x_j$.
\begin{remark}\label{stbil}If the action of $\,\Gm\hs$ on $\,\hm\,$ is free 
and properly dis\-con\-tin\-u\-ous, $\,a_j$ and $\,y_j$ are sequences in 
$\,\Gm\hs$ and $\,\hm$, while both $\,y_j$ and $\,a_jy_j$ converge, then the 
sequence $\,a_j$ is constant except for finitely many $\,j$.
\end{remark}
By a {\it lattice\/} in a real vector space $\,\lz\,$ with 
$\,\dim\lz<\infty\,$ we mean, as usual, an additive subgroup of $\,\lz\,$ 
generated by some basis of $\,\lz$.
\begin{remark}\label{lttce}For $\,\lz\,$ as above, a countable additive 
subgroup $\,\Lm\subset\lz$ is a lattice if and only if 
$\,\hs\text{\rm span}\,\Lm=\lz\,$ and $\,\Lm\,$ is closed as a 
subset of $\,\lz$. See \cite[Chap.\ VII, Th\'eor\`eme 2]{bourbaki}.
\end{remark}
\begin{lemma}\label{polyn}If\/ $\,k,l\in\bbZ\,$ and\/ 
$\,2\le k<l\le k^2\nnh/4$, then the polynomial 
$\,P(\lambda)=-\lambda^3\nh+k\lambda^2\nh-l\lambda+1\,$ in the real variable 
$\,\lambda\,$ has three distinct real roots $\,\lambda,\mu,\nu\,$ such that\/ 
$\,1/l\,<\,\lambda\,<\,1\,<\mu\,<\,k/2\,<\,\nu\,<\,k$, and hence
\begin{equation}\label{lmn}
0\,<\,\lambda\,<\,\mu\,<\,\nu\hs,\hskip13pt\lambda\,<\,1\,<\,\nu\hs,\hskip13pt
\lambda\mu\,<\,1\,<\,\mu\nu\hs,\hskip13pt\lambda\nu\,\ne\,1\hs.
\end{equation}
\end{lemma}
\begin{proof}Since $\,P(\lambda)=(k-\lambda)\lambda^2\nh+1-l\lambda$, we get 
$\,P(\lambda)\ge1-l\lambda>0$ if $\,\lambda<1/l\,$ and 
$\,P(\lambda)\le1-l\lambda\le1-kl<0\,$ if $\,\lambda\ge k$, as well as 
$\,P(\lambda)=(k-1/l)\lambda^2\nh>0\,$ for $\,\lambda=1/l$. Similarly, 
$\,P(1)=k-l<0\,$ and $\,P(\lambda)=1+(\lambda^2\nh-l)\lambda\ge1\,$ if 
$\,\lambda=k/2$. Therefore, $\,P\,$ has some roots $\,\lambda,\mu,\nu\,$ such 
that $\,1/l<\lambda<1<\mu<k/2<\nu<k$. This yields (\ref{lmn}): the 
inequalities $\,\lambda\mu<1\ne\lambda\nu\,$ follow since 
$\,\lambda\mu\nu=1\,$ and $\,\nu>1\ne\mu$.
\end{proof}
It is obvious that, for a nonzero polynomial $\,P\,$ in the real variable $\,\lambda$,
\begin{equation}\label{pol}
\begin{array}{l}
\text{\rm all\ coefficients\ of\ $\,P\hs$\ are\ integers,\ its\ leading\ term\ 
is\ $\,(-\lambda)^d\nh$,}\\
\text{\rm where\hskip.8pt\ $\,\hs d\,=\,\hs\text{\rm deg}\,P$,\hskip.8pt\ 
and\hskip.8pt\ its\hskip.8pt\ constant\hskip.8pt\ term\hskip.8pt\ 
equals\hskip.8pt\ $\,1\,$\hskip.8pt\ or\ $\,-1$}
\end{array}
\end{equation}
if and only if (\ref{pol}) holds with $\,P\,$ replaced by the product 
$\,(1-\lambda)\hs P$.

Let $\,T:\lz\to\lz\,$ be an endomorphism in a \rvs\ $\,\lz\,$ with 
$\,\dim\lz=d<\infty$. Clearly, $\,T(\Lm)=\Lm\,$ for some lattice $\,\Lm\,$ in 
$\,\lz\,$ if and only if $\,\text{\rm det}\,T=\pm1\,$ and the matrix of 
$\,T\,$ in some basis of $\,\lz\,$ consists of integers. Condition (\ref{pol}) 
characterizes the characteristic polynomials of those endomorphisms $\,T\,$ of 
$\,\lz\,$ for which such a lattice $\,\Lm\,$ exists.

In fact, given $\,a_1,\dots,a_d\in\bbR\hs$, let $\,\mathfrak{T}\,$ be the 
$\,d\times d\,$ matrix with the rows 
$\,[\hs a_1\hskip6pt\ldots\hskip6pta_d\hh]$, 
$[\hs1\hskip6pt0\hskip6pt\ldots\hskip6pt0\hh]$, 
$[\hs 0\hskip6pt1\hskip6pt0\hskip6pt\ldots\hskip6pt0\hh]$, $\ldots\,$, 
$[\hs 0\hskip6pt\ldots\hskip6pt0\hskip6pt1\hskip6pt0\hh]$. Then 
$\,\mathfrak{T}$ has the characteristic polynomial 
$\,(-1)^d(\lambda^d\nh-a_1\lambda^{d-1}\nnh-\ldots-a_{d-1}\lambda-a_d)$, 
which proves sufficiency of (\ref{pol}).

\section{Conformally symmetric Ric\-ci-re\-cur\-rent metrics}\label{csrr}
In this section, $\,\fh,p,n,\mv\nh,\lr\,$ and $\,A\,$ stand for the 
following objects:
\begin{equation}\label{obj}
\begin{array}{l}
\text{\rm a\ nonconstant\ periodic\ function\ 
$\fh\nnh:\nh\bbR\tci\bbR\hs$\ with\ a\ period\ $\hs p>0$,}\hskip-9pt\\
\text{\rm an\ integer\ $\,n\ge4\,$\ and\ a\ \rvs\ $\,\mv\hs$\ 
of\ dimension\ $\,n-2\hh$,}\hskip-9pt\\
\text{\rm a\ pseu\-\hbox{do\hs-}Euclid\-e\-an\ inner\ 
product\ $\,\lr\,$\ on\ $\,\mv\nnh$,}\\
\text{\rm a\ nonzero,\ traceless,\ $\,\lr$-self-ad\-joint\ linear\ operator\ 
$\,A:\mv\to\mv\nh$.\hskip-9pt}
\end{array}
\end{equation}
As in \cite{roter}, the data (\ref{obj}) lead to a \prc
\begin{equation}\label{met}
\hskip-10pt\hg\,=\,\kx\,dt^2\hs+\,dt\,ds\,+\,\vh\hskip14pt\text{\rm on\ the\ 
manifold\ $\,\hm=\rto\nh\times\mv\hs\approx\,\rn\nnh$.}
\end{equation}
The products of differentials stand here for symmetric products, $\,t,s\,$ are 
the Cartesian coordinates on $\,\rto$ treated, with the aid of the projection 
$\,\hm\to\rto\nnh$, as functions $\,\hm\to\bbR\hs$, and $\,\vh\,$ is the 
pull\-back to $\,\hm\,$ of the flat (constant) \prc\ on $\,\mv$ formed by the 
inner product $\,\lr$, while $\,\kx:\hm\to\bbR\hs$, with 
$\,\kx(t,s,v)=\fh(t)\hh\lg v,v\rg+\lg Av,v\rg$.
\begin{lemma}\label{ecsrr}For any choice of the data {\rm(\ref{obj})}, the 
metric $\,\hg\,$ given by {\rm(\ref{met})} is \ecs\ and \rr.
\end{lemma}
\begin{proof}See \cite[Theorem 3]{roter}, where weaker assumptions are used: 
rather than being defined on $\,\bbR\,$ and periodic, $\,\fh\,$ is just 
a nonconstant real-val\-ued $\,C^\infty$ function on an open interval 
$\,I\subset\bbR\hs$, and $\,\hm\,$ is replaced with 
$\,I\times\bbR\times\mv\nh$. In the notation of \cite{roter}, our $\,t,s\,$ 
and $\,\vh\,$ appear as $\,x^1\nh,\hs2x^n$ and 
$\,k_{\lambda\mu}\hs dx^\lambda\hs dx^\mu\nh$, while our $\,\fh(t)\,$ is 
$\,[2(n-2)]^{-1}\hs C\exp\,(\hs\int Q\,dx^1)$.
\end{proof}
Next, we set $\,\gp=\bbZ\times\bbR\times\xe$, where $\,\xe\hs$ is the vector 
space of all $\,C^\infty$ solutions $\,u:\bbR\to\mv\,$ to the differential 
equation $\,\ddot u(t)=\fh(t)\hh u(t)+Au(t)$. Clearly, 
$\,\varOmega(u,w)=\lg\dot u,w\rg-\lg u,\dot w\rg\,$ is, for any 
$\,u,w\in\xe$, a constant function $\,\bbR\to\bbR\hs$, which defines a 
nondegenerate \skc\ bi\-lin\-e\-ar form $\,\varOmega:\xe\times\xe\to\bbR\hs$. 
On the other hand, setting $\,(Tu)(t)=u(t-p)$, we obtain a linear isomorphism 
$\,T:\xe\to\xe\,$ with
\begin{equation}\label{tao}
\text{\rm $T\hh^*\nnh\varOmega=\varOmega$,\hskip13ptthat\ is,
\hskip6pt$\,\varOmega(Tu,Tw)=\varOmega(u,w)\,$\ whenever\ $\,u,w\in\xe$.}
\end{equation}
For $\,(k,q,u),(l,r,w)\in\gp\,$ and $\,(t,s,v)\in\hm=\rto\nh\times\mv\nh$, 
we set
\begin{equation}\label{qku}
\begin{array}{rl}
\mathrm{a)}&\hskip-2pt(k,q,u)\cdot(l,r,w)
=(k+l,\hs q+r-\varOmega(u,T^{\hs l\nh}w),\hs T^{\hs-\hs l\nh}u+w)\hs,\\
\mathrm{b)}&\hskip-2pt(k,q,u)\cdot(t,s,v)
=(t+kp,\hs s+q-\lg\dot u(t),2v+u(t)\rg,\hs v+u(t))\hs,\hskip-11pt
\end{array}
\end{equation}
which, by (\ref{tao}), defines a Lie-group structure in $\,\gp\,$ and an 
action of the Lie group $\,\gp\,$ on the manifold $\,\hm$. With all triples 
assumed to be elements of $\,\gp$, (\ref{qku}.a) gives
\begin{equation}\label{spc}
\begin{array}{rl}
\mathrm{i)}&\hskip0pt(k,q,u)^{-1}\nh=(-\hs k,-\hs q,-\hs T^{\hs k}u)\hs,\\
\mathrm{ii)}&\hskip0pt(k,q,u)\cdot(0,r,0)
=(0,r,0)\cdot(k,q,u)\,=\,(k,\,q+r,\,u)\hs,\\
\mathrm{iii)}&\hskip0pt(0,r,0)^l\nh\cdot(k,q,u)=(k,q+lr,u)\hs,\\
\mathrm{iv)}&\hskip0pt(k,q,u)\cdot(0,r,w)\cdot(k,q,u)^{-1}\nh
=(0,r-2\hh\varOmega(u,w),T^{\hs k}w)\hs,\\
\mathrm{v)}&\hskip0pt(0,q,u)\cdot(0,r,w)\cdot(0,q,u)^{-1}\nnh\cdot
(0,r,w)^{-1}\nh=(0,2\hh\varOmega(w,u),0)\hs.\hskip-8pt
\end{array}
\end{equation}
Our $\,\gp\,$ also acts on the manifold $\,\rto\nh\times\xe\hs$, \feic\ 
to $\,{\mathbf{R}}^{\hskip-.5pt2n-2}\nnh$, by
\begin{equation}\label{qkx}
(k,q,u)\cdot(t,z,w)\,=\,(t+kp,\,z+q-\varOmega(u,w),\,T^k(w+u))\hs.
\end{equation}
The following mapping is easily verified to be equivariant relative to the 
actions of $\,\gp\,$ given by (\ref{qkx}) and (\ref{qku}.b):
\begin{equation}\label{mpg}
\rto\nh\times\xe\ni(t,z,w)\mapsto(t,s,v)
=(t,\,z-\lg\dot w(t),w(t)\rg,\,w(t))\in\hm\hs.
\end{equation}
\begin{lemma}\label{isome}The group\/ $\,\gp\hs$ acts on $\,(\hm,\hg\hs)\,$ 
by isometries.
\end{lemma}
\begin{proof}Using any fixed basis $\,e_\lambda$ of $\,\mv\nh$, 
$\,\lambda=3,\dots,n$, we obtain $\,\vh=\vh_{\lambda\mu}\,dv^\lambda\hs dv^\mu$ 
and $\,\kx=(\fh \vh_{\lambda\mu}+a_{\lambda\mu})\hs v^\lambda v^\mu$ in the 
product coordinates $\,t,s,v^\lambda$ for $\,\hm$, where the coordinate 
functions $\,v^\lambda$ on $\,\mv$ send each $\,v\in\mv$ to its components in 
the expansion $\,v=v^\lambda e_\lambda$, while 
$\,\vh_{\lambda\mu}=\lg e_\lambda,e_\mu\rg\,$ and 
$\,a_{\lambda\mu}=\lg Ae_\lambda,e_\mu\rg$. For any given $\,(k,q,u)\in\gp$, 
the mapping $\,F:\hm\to\hm$ with $\,F(t,s,v)=(k,q,u)\cdot(t,s,v)\,$ has the 
components $\,F^*\nh t,F^*\nnh s,F^*\nh v^\lambda$ (that is, $\,t\circ F$, 
etc.) equal to $\,t+kp$, 
$\,s+q-\vh_{\lambda\mu}\dot u^\lambda(t)[2v^\mu\nh+u^\mu(t)]$ and 
$\,F^\lambda\nh=v^\lambda\nh+u^\lambda(t)$. Evaluating their 
differentials and noting that $\,\ddot u(t)=\fh(t)\hh u(t)+Au(t)$, we get 
$\,F^*\hg=(F^*\nh\kx)\hs(dF^*\nh t)^2+\hs(dF^*\nh t)\,dF^*\nnh s
+\vh_{\lambda\mu}\,dF^\lambda\,dF^\mu=\hg$, as required.
\end{proof}
\begin{remark}\label{trnsl}We will use the symbol $\,T\,$ for a more general 
translation operator, acting on functions $\,\bbR\ni t\mapsto\eta(t)\,$ valued 
in scalars, vectors or operators by $\,(T\eta)(t)=\eta(t-p)$.
\end{remark}

\section{First-or\-der and Lagrangian subspaces}\label{fols}
We again assume $\,\fh,p,n,\mv\nh,\lr\,$ and $\,A\,$ to be as in (\ref{obj}), 
while $\,\hm,\Gm$, $\,\xe,\hh\varOmega\,$ and $\,T\,$ stand for the 
corresponding objects defined in Section~\ref{csrr}.

By a {\it first-or\-der subspace} of the solution space $\,\xe\hs$ we mean 
any $\,(n-2)$\diml\ vector subspace $\,\lz\subset\xe\,$ having the 
property that $\,u(t)\ne0$ whenever $\,u\in\lz\smallsetminus\{0\}\,$ and 
$\,t\in\bbR\hs$. For any first-or\-der subspace $\,\lz$,
\begin{equation}\label{evl}
\begin{array}{l}
\text{\rm the\hs\hs\ evaluation\hs\hs\ operators\hs\ $\,\,u\mapsto u(t)\,\,$\ 
form\hs\hs\ a\hs\hs\ $\,\hs C^\infty$\hs\hs\ curve,}\\
\text{\rm parametrized\ by\ $\,t\in\bbR\hs$,\ of\ linear\ isomorphisms\ 
$\,\lz\to\mv\nnh$.}
\end{array}
\end{equation}
Since $\,\varOmega\,$ is nondegenerate, $\,\dim\lz'\nnh=\dim\hs\xe-\hh\dim\lz\,$ 
for any vector subspace $\,\lz\subset\xe\,$ and 
$\,\lz'=\{u\in\xe:\varOmega(u,w)=0\hskip4pt\text{\rm for\ all\ }\,w\in\lz\}$. 
Thus
\begin{equation}\label{dml}
\text{\rm $2\hs\dim\lz\,\le\,\dim\hs\xe\,\,$\ whenever\ $\,\lz\,$\ is\ 
a\ Lagrangian\ subspace\ of\ $\,\xe$,}
\end{equation}
where $\,\lz\,$ is called {\it Lagrangian\/} if $\,\varOmega(u,w)=0\,$ for all 
$\,u,w\in\lz$.
\begin{remark}\label{equiv}If $\,\lz\subset\xe\,$ is a first-or\-der subspace, 
the restriction of the mapping (\ref{mpg}) to $\,\rto\nh\times\lz\,$ is a 
\feo\ $\,\rto\nh\times\lz\to\hm$, equivariant relative to the actions 
(\ref{qkx}) and (\ref{qku}.b) of the subgroup $\,\hp\,$ of $\,\gp\,$ whose 
underlying set is $\,\{0\}\times\bbR\times\lz$. In fact, by (\ref{evl}), the 
restriction is a \feo, and it is equivariant since so is (\ref{mpg}).
\end{remark}
\begin{lemma}\label{rgsbs}First-or\-der subspaces\/ $\,\lz\hs$ of the solution 
space\/ $\,\xe$ are in a bijective correspondence with $\,C^\infty\nnh$ 
functions $\,B:\bbR\to\hs\text{\rm End}\hs(\mv)\,$ such that\/ 
$\,\dot B+B^2=\fh+A$, where $\,\fh\,$ stands for the function 
$\,t\mapsto\fh(t)\,\text{\rm Id}\hh$. The correspondence assigns to\/ $B\hs$ 
the space\/ $\,\lz\hs$ of all solutions\/ $\,u:\bbR\to\mv\hs$ to the 
differential equation $\,\dot u(t)=B(t)\hh u(t)$, and, for this $\,\lz$,
\begin{enumerate}
  \def\theenumi{{\rm\roman{enumi}}}
\item $\lz\,$ is a Lagrangian subspace of\/ $\,\xe\hs$ if and only if\/ 
$\,B(t)\,$ is self-ad\-joint relative to $\,\lr\,$ for every $\,t\in\bbR\hs$,
\item $\lz\,$ is $\,T\nnh$\inv\ if and only if\/ $\,B\,$ is periodic with 
period\/ $\,p$,
\item if\/ $\hs T(\lz)\nh=\nh\lz$, the determinant of\/ 
$\hs T\hskip-2pt:\hskip-1pt\lz\to\lz\,$ is\/ 
$\hs\exp\,(\nh-\hskip-3.3pt\int_0^p\text{\rm tr}\hskip1ptB(t)\,dt)$.
\end{enumerate}
\end{lemma}
\begin{proof}The assignment $\,B\mapsto\lz\,$ described in the lemma sends 
$\,B\,$ with $\,\dot B+B^2=\fh+A\,$ to a first-or\-der subspace of $\,\xe\hs$ 
in view of uniqueness of solutions for ordinary differential equations. The 
surjectivity and injectivity of $\,B\mapsto\lz\,$ are both obvious from 
(\ref{evl}): given $\,\lz$, we choose $\,B(t)\,$ to be the inverse of the 
evaluation isomorphism $\,u\mapsto u(t)\,$ followed by the operator 
$\,u\mapsto\dot u(t)$, which is clearly the unique choice of $\,B$ producing 
the given solution space $\,\lz$. Now (i) is immediate from (\ref{evl}).

The equation $\,\dot u=B\hh u\,$ for $\,u:\bbR\to\mv\hs$ implies 
$\,(Tu)\dot{\,}=(TB)\hh(Tu)\,$ (cf.\ Remark~\ref{trnsl}). This second equation 
gives $\,T(\lz)\subset\lz\,$ whenever $\,B$ is periodic with period $\,p\,$ 
(that is, $\,TB=B$). Conversely, if $\,T(\lz)=\lz$, the two equations 
combined with (\ref{evl}) yield $\,TB=B$, proving (ii).

For a basis $\,u_\lambda\,$ of $\,\lz$, $\,\lambda=3,\dots,n$, and a fixed 
volume form $\,[\,\ldots\,]\,$ in $\,\mv\nh$, defining $\,\eta:\bbR\to\bbR\,$ 
by $\,\eta=[u_3,\ldots,u_n]\,$ we get 
$\,\dot\eta=\eta\,\hs\text{\rm tr}\,B$, so that, if $\,\lz\,$ is 
$\,T\nnh$\inv, integration shows that $\,\log|\eta|-\log|T\eta|=\tau$, where 
$\,\tau=\int_0^p\text{\rm tr}\,B(t)\,dt\,$ (cf.\ (ii) and 
Remark~\ref{trnsl}), and, therefore, $\,\eta=e^\tau T\eta$. On the other 
hand, $\,\eta\,\text{\rm det}\,T$, for $\,T:\lz\to\lz$, equals 
$\,[Tu_3,\ldots,Tu_n]=T\eta$. Hence 
$\,\text{\rm det}\,T=\eta^{-1}T\eta=e^{-\hs\tau}\nnh$, which gives (iii).
\end{proof}
\begin{remark}\label{cmmut}If $\,\lz\hs$ and $\,B\,$ are related as in 
Lemma~\ref{rgsbs}, while $\,T(\lz)=\lz\,$ and $\,B(t)\,$ commutes with 
$\,B(t\hh'\hh)$ for all $\,t,t\hh'\nh\in\bbR\hs$, then $\,T:\lz\to\lz$ is 
given by $\,e^{-\hs S}\nnh$, where 
$\,S=\int_0^p B(t)\,dt\in\hs\text{\rm End}\hs(\mv)\,$ acts as an operator 
$\,\lz\to\lz\,$ with $\,(Su)(t)=Su(t)$.

In fact, let $\,J(t,s)=\int_s^tB(t\hh')\,dt\hh'\in\hs\text{\rm End}\hs(\mv)\,$ 
for $\,t,s\in\bbR\hs$. Since $\,d\hs e^{J(t,s)}\nh/dt=B(t)\hs e^{J(t,s)}\nnh$, 
the unique solution $\,w\in\lz\,$ to the initial value problem $\,\dot w=Bw$, 
$\,w(s)=v\,$ is, for any $\,s\in\bbR\,$ and $\,v\in\mv\nh$, given by 
$\,w(t)=e^{J(t,s)}v$. Applying this to $\,w=u\,$ or $\,w=Tu$, where 
$\,u\in\lz$ is fixed, and $\,s=0$, we see that $\,u(t)=e^{J(t,0)}u(0)\,$ 
and $\,(Tu)(t)=e^{J(t,p)}u(0)$, for each $\,t\in\bbR\,$ (since $\,(Tu)(p)=u(0)$). 
Thus, $\,(Tu)(t)=e^{J(t,p)-J(t,0)}u(t)=e^{-J(0,p)}u(t)$.
\end{remark}

\section{Discrete subgroups of $\,\gp$}\label{dsog}
For $\,\fh,p,n,\mv\nh,\lr,A\,$ and $\,\xe,\hh\varOmega,T,\gp\,$ as in 
Section~\ref{csrr}, let
\begin{equation}\label{pgz}
\varPi:\gp\to\bbZ\hs,\hskip8pt\Dt:\text{\rm Ker}\,\varPi\to\xe\hs,\hskip6pt
\mathrm{and}\hskip6pt\Dt_t:\text{\rm Ker}\,\varPi\to\mv\nnh,\hskip5pt
\mathrm{for\ }t\in\bbR\hs,
\end{equation}
be the following homomorphisms (with $\,\xe,\mv\hs$ treated as additive groups):
\begin{equation}\label{hom}
\varPi(k,q,u)=k\hs,\hskip15pt\Dt(0,q,u)=u\hs,\hskip15pt
\Dt_t(0,q,u)=u(t)\hs.
\end{equation}
Throughout this section, given a subgroup $\,\Gm\hs$ of $\,\gp$, we use the 
notation
\begin{equation}\label{llx}
\gkp\,=\,\Gm\hskip1pt\cap\hskip1pt\text{\rm Ker}\,\varPi\hs,\hskip8pt
\Xi\,=\,\gkp\hh\cap\text{\rm Ker}\,\Dt\hs,\hskip8pt
\Lm\,=\,\Dt(\gkp)\hs,\hskip8pt
\lz\,=\,\hs\text{\rm span}\,\Lm\hs.
\end{equation}
Thus, $\,\Xi\,$ and $\,\Lm\,$ are the kernel and image of $\,\Dt:\gkp\to\xe$, 
while $\,\lz\,$ is the vector subspace of $\,\xe\hs$ spanned by the additive 
subgroup $\,\Lm$. Since 
$\,\,\gkp\,\subset\,\{0\}\times\bbR\times\Lm\,
\subset\,\{0\}\times\bbR\times\lz$, identifying 
$\,\{0\}\times\bbR\times\lz\,$ with $\,\bbR\times\lz$
\begin{equation}\label{tre}
\text{\rm we\ treat\ $\,\gkp\,$\ as\ a\ subset\ of\ $\,\bbR\times\lz\,$\ such\ 
that\ $\,\gkp\subset\bbR\times\Lm$.}
\end{equation}
\begin{remark}\label{ntrst}If a subgroup $\,\Gm\subset\gp\,$ acts on 
$\,\hm=\rto\nh\times\mv\hs$ freely and properly dis\-con\-tin\-u\-ous\-ly with 
a compact quotient\/ $\,M=\hm\nnh/\hh\Gm$, then
\begin{enumerate}
  \def\theenumi{{\rm\roman{enumi}}}
\item the image $\,\varPi(\Gm)\,$ equals $\,m\bbZ\,$ for some integer 
$\,m>0$,
\item $M\,$ is the total space of a $\,C^\infty$ bundle over the circle 
$\,\bbR/\nnh\varPi(\Gm)$, and the mapping 
$\,\hm\ni(t,s,v)\mapsto p^{-1}t\in\bbR\,$ descends to the bundle projection 
$\,\hs\text{\rm pr}:\hm\to\bbR/\nnh\varPi(\Gm)$,
\item for every $\,t\in\bbR\hs$, the submanifold 
$\,\hmt=\{t\}\times\bbR\times\mv\hs$ of $\,\hm\,$ is invariant under the 
subgroup $\,\gkp\,$ of $\,\Gm$, the action of $\,\gkp\,$ on $\,\hmt$ is 
properly dis\-con\-tin\-u\-ous, and the inclusion $\,\hmt\to\hm\,$ descends to 
an embedding $\,\hmt/\hh\gkp\to M$, the image of which is the fibre 
$\,M_{\text{\rm pr}\hh(t)}$ of the bundle $\,M\,$ over the point 
$\,\hs\text{\rm pr}\hh(t)\,$ in the base circle,
\item $\gkp\,$ acts on each $\,\hmt$ by affine transformations whose linear 
parts preserve the vector $\,(0,1,0)\in\{0\}\times\bbR\times\mv\subset\hm\,$ 
(cf.\ (\ref{qku}.b)); this gives rise to a flat tor\-sion\-free connection 
with a nonzero parallel vector field on the compact fibre 
$\,M_{\text{\rm pr}\hh(t)}$.
\end{enumerate}
In fact, (i) and (ii) are obvious as the assignment $\,(t,s,v)\mapsto t/p\,$ 
descends to a surjective submersion $\,M\to\bbR/\nnh\varPi(\Gm)$, which 
would be an unbounded $\,C^\infty$ function $\,M\to\bbR\,$ if 
$\,\varPi(\Gm)\,$ were the trivial group. As (\ref{qku}.b) clearly implies 
$\,\gkp$-in\-var\-i\-ance of $\,\hmt$, assertion (iii) is immediate from the 
definition of proper dis\-con\-ti\-nu\-i\-ty in Section~\ref{prel}.
\end{remark}
\begin{theorem}\label{conse}Suppose that\/ $\,\Gm\hs$ is a subgroup of\/ 
$\,\gp$, the action of\/ $\,\Gm\hs$ on $\,\hm\,$ is free and properly 
dis\-con\-tin\-u\-ous, and the quotient manifold\/ $\,M=\hm\nnh/\hh\Gm\,$ is 
compact. Then
\begin{enumerate}
  \def\theenumi{{\rm\alph{enumi}}}
\item for every $\,t\in\bbR\hs$, the image $\,\Dt_t(\gkp)\,$ spans 
$\,\mv$ as a vector space,
\item $\gkp\hh\cap\text{\rm Ker}\,\Dt_t\hs=\,\Xi\,\,$ whenever 
$\,\,t\in\bbR\hs$,
\item $\Xi=\{0\}\times\bbZ\rc\times\{0\}=\{(0,l\rc,0)\nnh:\hs l\in\bbZ\}\,$ 
for some\/ $\,\rc\in[\hs0,\infty)$, 
\item $2\hh\varOmega(u,w)\in\bbZ\rc\,$ for\/ $\,\rc\,$ defined in\/ 
{\rm(c)} and all\/ $\,u,w\in\Lm$,
\item whenever $\,(k,q,u)\,$ is an element of\/ $\,\Gm$, we have 
$\,T^{\hh k}(\Lm)=\Lm\,$ and\/ $\,T^{\hh k}(\lz)=\lz$, while 
$\,\varPsi(\gkp)=\gkp$, for $\,\gkp\subset\bbR\times\lz\,$ as in 
{\rm(\ref{tre})} and\/ $\,\varPsi:\bbR\times\lz\to\bbR\times\lz\,$ given 
by $\,\varPsi(r,w)=(r-2\hh\varOmega(u,w),T^{\hs k}w)$,
\item $\Gm\,$ has no Abelian subgroup of finite index, unless 
$\,T^{\hs k}:\lz\to\lz$ equals the identity for some $\,(k,q,u)\in\Gm\,$ 
with $\,k\ge1$,
\item $\lz\,$ is a first-or\-der subspace of\/ $\,\xe$, so 
that\/ $\,\dim\lz=n-2$,
\item one of the following two cases occurs\/{\rm:}
\begin{enumerate}
  \def\theenumi{{\rm\alph{enumi}}}
\item[I)] $\Xi\,$ is the trivial group, $\,\lz\subset\xe\,$ is a Lagrangian 
subspace, $\,\gkp$ with {\rm(\ref{tre})} is a lattice in 
$\,\bbR\times\lz$, and\/ $\,\Lm\,$ is the isomorphic image of\/ $\,\gkp\,$ 
under the projection $\,\bbR\times\lz\to\lz$,
\item[II)] $\Xi\,$ is isomorphic to $\,\bbZ\,$ and\/ $\,\Lm\,$ is a 
lattice in $\,\lz$.
\end{enumerate}
\end{enumerate}
\end{theorem}
\begin{proof}If $\,\Dt_t(\gkp)\,$ spanned a proper subspace 
$\,\mv\hh'$ of $\,\mv\nh$, a nonzero linear functional $\,\mv\nh\to\bbR\,$ 
vanishing on $\,\mv\hh'$ would descend to an unbounded $\,C^\infty$ function 
on the compact manifold $\,M_t$ (cf.\ Remark~\ref{ntrst}(iii)). This yields 
(a). Next, if we had $\,u(t)=0\,$ for some $\,(0,q,u)\in\gkp\,$ with 
$\,u\in\lz\smallsetminus\{0\}\,$ and some $\,t\in\bbR\hs$, choosing 
$\,v\in\mv$ such that $\,q=\lg\dot u(t),2v\rg$ (which exists as 
$\,\dot u(t)\ne0$), we would get $\hs(0,q,u)\cdot(t,s,v)=(t,s,v)\hs$ with any 
$\,s\in\bbR\,$ (cf.\ (\ref{qku}.b)), and so the action of $\,\Gm\hs$ would 
not be free. Hence, if $\,(0,q,u)\in\gkp\hh\cap\text{\rm Ker}\,\Dt_t$, then 
$\,u=0$, which implies (b).

Obviously, $\,\Xi\hs=\{0\}\nh\times\kp\times\{0\}\,$ for some additive 
subgroup $\,\kp\hh$ of $\,\bbR\hs$. If $\,\hs\kp$ were not a closed subset of 
$\,\bbR\hs$, any fixed sequence $\,\rc_j$ of mutually distinct nonzero 
elements of $\,\hs\kp$ with $\,\rc_j\to0$, and any $\,(t,s,v)\in\hm$, would 
yield the sequences $\,(0,\rc_j,0)\,$ and $\,(0,\rc_j,0)\cdot(t,s,v)\,$ 
contradicting, by (\ref{spc}.ii), the conclusion of Remark~\ref{stbil}. This 
proves (c) (cf.\ Remark~\ref{lttce}).

By (\ref{spc}.v), the commutator of $\,(0,q,u),(0,r,w)\in\gkp\,$ is 
$\,(0,2\hh\varOmega(w,u),0)$, which, in view of (\ref{llx}), is an element of 
$\,\Xi$. Thus, (c) gives (d).

Next, let $\,(k,q,u)\in\Gm$. By (\ref{spc}.iv) and (\ref{spc}.i), the inner 
automorphisms corresponding to $\,(k,q,u)\,$ and $\,(k,q,u)^{-1}$ send any 
$\,(0,r,w)\in\gkp\,$ to $\,(0,r-2\hh\varOmega(u,w),T^{\hs k}w)\,$ and 
$\,(0,r+2\hh\varOmega(T^{\hs k}u,w),T^{\hs-\hs k}w)$, which must again be 
elements of the normal subgroup $\,\gkp\subset\Gm$. Thus, both $\,T^{\hs k}$ 
and $\,T^{\hs-\hs k}$ leave $\,\Lm\,$ invariant, which yields (e), as 
$\,\lz=\hs\text{\rm span}\,\Lm\,$ and 
$\,\varOmega(T^{\hs k}u,w)=\varOmega(u,T^{\hs-\hs k}w)\,$ by (\ref{llx}) and 
(\ref{tao}).

If $\,\Gm\hs$ has an Abelian subgroup $\,\Gm\hs'$ of finite index, replacing 
$\,\Gm\hs$ by $\,\Gm\hs'\nnh$, we may assume that $\,\Gm\hs$ is Abelian. From 
(\ref{spc}.iv) for any $\,(0,r,w)\in\gkp$ and any $\,(k,q,u)\in\Gm\,$ 
with $\,k\ge1\,$ (which exists, cf.\ Remark~\ref{ntrst}(i)) we now get, by 
(\ref{llx}), $\,T^{\hs k}w=w\,$ for all $\,w\in\lz=\hs\text{\rm span}\,\Lm$, 
and (f) follows.

The evaluation operator $\,\lz\to\mv\hs$ is surjective for each $\,t$, since 
it maps $\,\Lm=\Dt(\gkp)\,$ onto 
$\,\Lm_t=\Dt_t(\gkp)$, while $\,\lz=\hs\text{\rm span}\,\Lm\,$ and 
$\,\mv\nh=\hs\text{\rm span}\,\Lm_t$ by (\ref{llx}) and (a). Consequently, 
$\,\dim\lz\ge\dim\hs\mv\nh=n-2$. Proving (g) is now reduced to 
showing that $\,\dim\lz\,\le\hs n-2$.

We first assume that $\,\Xi\,$ is trivial, and so $\,\Dt:\gkp\to\Lm\,$ 
is an isomorphism. Thus, $\,\gkp\,$ is Abelian, and (\ref{qku}.a) with 
$\,k=l=0\,$ implies that $\,\lz=\hs\text{\rm span}\,\Lm\,$ is a Lagrangian 
subspace, which has two consequences. First, (g) holds in this case, as 
$\,\dim\lz\le n-2\,$ by (\ref{dml}). Secondly, again by (\ref{qku}.a), the 
subgroup $\,\hp\,$ of $\,\gp\,$ with the underlying set 
$\,\{0\}\times\bbR\times\lz\,$ is Abelian, and, under the obvious 
identification $\,\{0\}\times\bbR\times\lz\,\approx\,\bbR\times\lz$, it 
coincides with the additive group of $\,\bbR\times\lz$. Hence $\,\gkp\,$ 
with (\ref{tre}) is an additive subgroup of $\,\bbR\times\lz$. The properly 
dis\-con\-tin\-u\-ous action of $\,\gkp\,$ on $\,\hmt$, having the compact 
quotient $\,M_t$ (see Remark~\ref{ntrst}(iii)) corresponds, under the 
equivariant \feo\ defined in Remark~\ref{equiv} (for our $\,\lz$) combined 
with the obvious identification 
$\,\{t\}\times\bbR\times\lz\,\approx\,\bbR\times\lz$, to the action of 
$\,\gkp\,$ by vec\-tor-space translations on the ambient space 
$\,\bbR\times\lz$. Therefore, $\,\gkp\,$ is a lattice in $\,\bbR\times\lz$, 
and assertion (h-I) follows.

Finally, let $\,\Xi\,$ be nontrivial. Thus, $\,\rc\,$ in (c) is positive. 
By (\ref{spc}.iii),
\begin{equation}\label{qlr}
\text{\rm $(k,q+l\rc,u)\in\Gm\,\,$\ whenever\ $\,(k,q,u)\in\Gm\,$\ and\ 
$\,l\in\bbZ$.}
\end{equation}
Consequently, $\,\Lm\,$ is closed as a 
subset of $\,\lz$. In fact, otherwise there would exist a sequence $\,u_j$ of 
mutually distinct nonzero elements of $\,\Lm$, such that $\,u_j\to0$. 
Choosing $\,(0,q_j,u_j)\in\gkp\,$ with $\,q_j\in[\hs0,\rc\hs]\,$ (cf.\ 
(\ref{qlr})), and then replacing the $\,q_j$ by a convergent subsequence, we 
could use any $\,(t,s,v)\in\hm\,$ to obtain a sequence 
$\,(0,q_j,u_j)\cdot(t,s,v)\,$ that converges (see (\ref{qku}.b)), contrary to 
Remark~\ref{stbil}. Similarly, $\,\Lm_t=\Dt_t(\gkp)$ is, for each $\,t$, a 
closed subset of $\,\mv\nh$. Namely, if for some $\,t\,$ and some sequence 
$\,(0,q_j,u_j)\in\gkp\,$ the vectors $\,u_j(t)\in\mv\hs$ were mutually 
distinct and nonzero, while $\,u_j(t)\to0$, (\ref{qlr}) would allow us to 
modify $\,q_j$ in such a way that $\,q_j-\lg\dot u_j(t),u_j(t)\rg\,$ is a 
bounded sequence in $\,\bbR\hs$. Any convergent subsequence of 
$\,(0,q_j,u_j)\cdot(t,s,0)$, with any $\,s\in\bbR\hs$, would now, by 
(\ref{qku}.b), again contradict Remark~\ref{stbil}.

Being closed, both $\,\Lm\,$ and $\,\Lm_t=\Dt_t(\gkp)\,$ are lattices in the 
respective spaces $\,\lz\,$ and $\,\mv\hs$ (see Remark~\ref{lttce} and (a)). 
Thus, they are free Abelian groups of ranks $\,\dim\lz\ge n-2\,$ and, 
respectively, $\,\dim\hs\mv\nh=n-2$. The restriction to $\,\Lm\,$ of the 
evaluation operator $\,\lz\to\mv\hs$ is, by  (\ref{hom}), a surjective 
homomorphism onto $\,\Lm_t$, and, by (b), it is injective. Hence 
$\,\dim\lz=\dim\hs\mv\nh$, so that we have (g),  while (c) gives (h-II). 
\end{proof}

\section{Simplifying assumptions}\label{smpl}
For any given objects (\ref{obj}), with $\,\hm\,$ and $\,\gp\,$ as in 
Section~\ref{csrr}, we ask whether some subgroup $\,\Gm\subset\gp\,$ acts on 
$\,\hm\,$ properly dis\-con\-tin\-u\-ous\-ly, producing a {\it compact\/} 
$\,n$\diml\ quotient manifold $\,M=\hm\nnh/\hh\Gm$. For the purpose of 
answering this question, we may always require $\,\Gm\hs$ to satisfy the 
following additional conditions:
\begin{equation}\label{smp}
\begin{array}{rl}
\mathrm{a)}&\hskip-3pt\text{\rm the\ integer\ $\,m\,$\ defined\ in\ 
Remark~\ref{ntrst}(i)\ is\ equal\ to\ $\,1$,}\\
\mathrm{b)}&\hskip-3pt\text{\rm $\varOmega(u,w)\in\bbZ\rc\,$\ for\ all\ 
$\,u,w\in\Lm\hs$\ (cf.\ (\ref{llx})\ and\ Theorem~\ref{conse}(c)),
\hskip-15pt}\\
\mathrm{c)}&\hskip-3pt\text{\rm $\gkp\,$\ is\ a\ lattice\ in\ 
$\,\bbR\times\lz$,\ where\ $\,\gkp\subset\bbR\times\lz\,$\ as\ in\ 
(\ref{tre}).}
\end{array}
\end{equation}
In fact, assuming that $\,m=1\,$ leads to no loss of generality, as the use of 
$\,mp\,$ instead of $\,p\,$ in (\ref{obj}) causes $\,m\,$ to be replaced by 
$\,1$.

Next, by (\ref{spc}.ii), the formula $\,\fy(r)=(0,r,0)\,$ defines a 
homomorphism $\,\fy\,$ from the additive group $\,\bbR\,$ into the center of 
$\,\gp$. The image $\,\fy\hs(\bbR)$ thus consists of isometries of 
$\,(\hm,\hg\hs)\,$ (cf.\ Lemma~\ref{isome}) descending to the quotient 
manifold $\,M=\hm\nnh/\hh\Gm$, so as to form a group $\,\kp\hs$ of isometries 
of $\,M$. The action of $\,\kp\hs$ on $\,M\,$ is free, and the kernel of the 
projection homomorphism $\,\fy\hs(\bbR)\to\kp\hs$ is the group $\,\Xi\,$ in 
(\ref{llx}). In 
fact, the kernel obviously contains $\,\Xi$. Now let $\,(t,s,v)\in\hm\,$ and 
$\,r\in\bbR\hs$. If $\,(0,r,0)\cdot(t,s,v)=(t,\hs s+r,\hs v)\,$ lies in the 
$\,\Gm$-or\-bit of $\,(t,s,v)$, and hence equals $\,(k,q,u)\cdot(t,s,v)\,$ for 
some $\,(k,q,u)\in\Gm$, then, by (\ref{qku}.b) and Theorem~\ref{conse}(b), 
$\,k=u=0\,$ and $\,q=r$, and so $\,(0,r,0)\in\hs\Xi$.

In case I of Theorem~\ref{conse}(h), (\ref{smp}.b) and (\ref{smp}.c) always 
hold, in view of Theorem~\ref{conse}(d) with $\,\rc=0$. To obtain 
(\ref{smp}.b) in case II of Theorem~\ref{conse}(h), we replace the quotient 
manifold $\,M=\hm\nnh/\hh\Gm\,$ by its own quotient $\,M\hh'$ under the free 
action of the $\,\bbZ_2$ subgroup of the circle $\,\kp=\fy\hs(\bbR)/\hs\Xi$. 
This means replacing $\,\Gm\hs$ by the subgroup 
$\,\Gm\hs'\nh=\Gm\cup(\zeta\nh\cdot\Gm)$ of $\,\gp\,$ generated by 
$\,\Gm\hs$ and the central element $\,\zeta=(0,\rc/2,0)$, for $\,\rc\,$ 
defined in Theorem~\ref{conse}(c). Since $\,\rc\hh'\nh=\rc/2\,$ corresponds to 
$\,\Gm\hs'$ just as $\,\rc\,$ did to $\,\Gm$, Theorem~\ref{conse}(d) now 
yields (\ref{smp}.b). Finally, (\ref{smp}.c) follows in case II from 
(\ref{smp}.b). In fact, as $\,\Lm=\Dt(\gkp)\,$ by (\ref{llx}) and $\,\Lm\,$ is 
a lattice in $\,\lz$ (see Theorem~\ref{conse}(h-II)), we may fix 
$\,(q_\lambda,u_\lambda)\in\gkp\subset\bbR\times\lz$, $\,\lambda=3,\dots,n$, 
such that $\,u_\lambda$ form a basis of $\,\lz\,$ generating the additive 
subgroup $\,\Lm$. We now show that $\,\gkp\,$ coincides with the lattice 
$\,\gkp\hs'$ in $\,\bbR\times\lz\,$ generated by the basis consisting of 
$\,(\rc\hh,0)\,$ and our $\,(q_\lambda,u_\lambda)$. Namely, the projection 
$\,(r,w)\mapsto w\,$ forms surjective group homomorphisms $\,\Dt:\gkp\to\Lm\,$ 
and $\,\gkp\hs'\nh\to\Lm$. Hence $\,\sum_\lambda k(\lambda)u_\lambda\,$ is 
the $\,\Dt$-im\-age of the product (in $\,\gkp$)
\begin{equation}\label{tol}
\text{\rm $(\rc\hh,0)^l(q_3,u_3)^{k(3)}\nnh\ldots\hs(q_n,u_n)^{k(n)}\in\gkp$,\ \ 
with\ integers\ $\,l\,$\ and\ $\,k(\lambda)$.}
\end{equation}
Thus, every $\,(r,w)\in\gkp\,$ is of 
the form (\ref{tol}): $\,(r,w)\,$ has the same $\,\Dt$-im\-age $\,w\,$ as 
$\,(\tilde r,w)=(q_3,u_3)^{k(3)}\nnh\ldots\hs(q_n,u_n)^{k(n)}$, with 
$\,k(\lambda)\in\bbZ\,$ such that $\,w=\sum_\lambda k(\lambda)u_\lambda$, 
and so $\,(r,w)=(\rc\hh,0)^l(\tilde r,w)\,$ for some 
$\,l\in\bbZ$. (By Theorem~\ref{conse}(c), $\,(\rc\hh,0)\,$ generates 
$\,\Xi$, the kernel of $\,\Dt:\gkp\to\Lm$.)

On $\,\bbR\times\Lm\,$ there are two group structures: one 
additive, with $\,\gkp\hs'$ as a subgroup, the other, given by 
$\,(q,u)(r,w)=(q+r-\varOmega(u,w),\hs u+w)$ (cf.\ (\ref{qku}.a)), 
having $\,\gkp\,$ as a subgroup. The mapping 
$\,\chi:\bbR\times\Lm\to S^1\nh\times\Lm$ with 
$\,\chi(r,w)=(r+\bbZ\rc,w)\,$ is, by (\ref{smp}.b), a homomorphism 
from {\it both} groups into the direct product of 
$\,S^1\nh=\bbR/\bbZ\rc\,$ and $\,\Lm$.

Any $\,(r,w)\in\gkp\,$ of the form (\ref{tol}) is related to the linear 
combination $\,(r\hs'\nnh,w)=l(\rc\hh,0)
+\sum_\lambda k(\lambda)(q_\lambda,u_\lambda)\in\gkp\hs'$ with the same 
coefficients $\,l,k(\lambda)\in\bbZ\,$ by 
$\,(r\hs'\nnh,w)=(r,w)+l\hh'\nh(\rc\hh,0)=(\rc\hh,0)^{l\hh'}(r,w)\,$ 
for some $\,l\hh'\nh\in\bbZ$. In fact, as 
$\,\chi(r,w)=\chi(r\hs'\nnh,w)\,$ (both being 
$\,\sum_\lambda k(\lambda)q_\lambda\,+\,\bbZ\rc$), $\,(r,w)\,$ and 
$\,(r\hs'\nnh,w)\,$ differ, relative to either group structure in 
$\,\bbR\times\Lm$, by an element of 
$\,\hs\text{\rm Ker}\hskip2pt\chi=\Xi=\bbZ\rc\times\{0\}$. 
Therefore, $\,\gkp=\gkp\hs'\nnh$.

\section{A criterion for the existence of compact quotients}\label{cecq}
Given $\,\fh,p,n,\mv\nh,\lr,A\,$ with (\ref{obj}), let 
$\,B:\bbR\to\hs\text{\rm End}\hs(\mv)\,$ be a $\,C^\infty$ function, periodic 
with period $\,p$, such that $\,\dot B+B^2=\fh+A$. These data lead to further 
objects: $\,\hm=\rto\nh\times\mv\hs$ and the group $\,\gp\,$ acting on 
$\,\hm\,$ defined in Section~\ref{csrr}, the vector space $\,\lz\,$ of 
dimension $\,n-2\,$ formed by all solutions $\,u:\bbR\to\mv\hs$ to the 
differential equation $\,\dot u(t)=B(t)\hh u(t)$, the translation operator 
$\,T:\lz\to\lz\,$ with $\,(Tw)(t)=w(t-p)$, the $\,(n-1)$\diml\ vector space 
$\,\wv=\bbR\times\lz$, and its one\diml\ subspace $\,\iv=\bbR\times\{0\}$.
\begin{theorem}\label{ncsuf}For any objects $\,\fh,p,n,\mv\nh,\lr\,$ and\/ 
$\,A\,$ as in {\rm(\ref{obj})}, the following two conditions are 
equivalent\/{\rm:}
\begin{enumerate}
  \def\theenumi{{\rm\roman{enumi}}}
\item some subgroup $\,\Gm\subset\gp\,$ acts on $\,\hm\,$ freely and properly 
dis\-con\-tin\-u\-ous\-ly with a compact quotient manifold\/ 
$\,M=\hm\nnh/\hh\Gm$, and satisfies condition {\rm(\ref{smp}.a);}
\item there exist a $\,C^\infty$ function 
$\,B:\bbR\to\hs\text{\rm End}\hs(\mv)$, periodic of period $\,p$, with 
$\,\dot B+B^2=\fh+A$, a lattice $\,\gkp\,$ in $\,\wv$, a linear functional\/ 
$\,\varphi\in\lz^*\nnh$, and\/ $\,\rc\in[\hs0,\infty)\,$ such that
\begin{enumerate}
  \def\theenumi{{\rm\alph{enumi}}}
\item[a)] $\gkp\cap\iv\,=\,\bbZ\rc\times\{0\}$,
\item[b)] $\varPsi(\gkp)=\gkp\,$ for $\,\varPsi:\wv\to\wv\,$ given by 
$\,\varPsi(r,w)=(r+\varphi(w),Tw)$,
\item[c)] $\varOmega(u,w)\in\bbZ\rc\,$ whenever $\,u,w\in\Lm$,
\end{enumerate}
where $\,\lz,T,\wv,\iv\,$ correspond to $\,B\,$ as above, $\,\Lm\,$ 
is the image of\/ $\,\gkp$ under the projection $\,\wv\to\lz$, and\/ 
$\,\varOmega(u,w)\in\bbR\,$ is the constant function 
$\,\lg\dot u,w\rg-\lg u,\dot w\rg$.
\end{enumerate}
\end{theorem}
We now show that (i) implies (ii) in Theorem~\ref{ncsuf}, postponing the 
proof of the converse statement until Section~\ref{ptio}.

Let $\,\Gm\hs$ be as in (i). According to Section~\ref{smpl}, we may also 
assume (\ref{smp}). By Theorem~\ref{conse}(e),\hs(g) with $\,k=1$, the 
space $\,\lz\,$ in (\ref{llx}) is a $\,T\nnh$\inv\ first-or\-der subspace 
of the solution space $\,\xe\,$ defined in Section~\ref{csrr}, and so, by 
Lemma~\ref{rgsbs}, $\,\lz\,$ arises from a $\,C^\infty$ function 
$\,B:\bbR\to\hs\text{\rm End}\hs(\mv)\,$ with $\,\dot B+B^2=\fh+A$. 
Lemma~\ref{rgsbs}(ii) shows that $\,B\,$ is periodic with period $\,p$. By 
(\ref{smp}.c), $\,\gkp\,$ defined in (\ref{llx}) is a lattice in 
$\,\wv=\bbR\times\lz$. Choosing $\,\rc\,$ as in Theorem~\ref{conse}(c), we 
obtain assertions (ii-a) and (ii-c), cf.\ (\ref{smp}.b). Finally, let us set 
$\,\varphi(w)=-\hs2\hh\varOmega(u,w)\,$ for $\,w\in\lz$, with $\,\varOmega\,$ 
defined in the lines preceding (\ref{qku}) and $\,u\in\xe\,$ chosen so that 
$\,(1,q,u)\in\Gm\,$ for some $\,q\in\bbR\hs$. (Note that we assume 
(\ref{smp}.a).) Now (ii-b) is immediate from Theorem~\ref{conse}(e) with 
$\,k=1$.
\begin{remark}\label{nilbd}As we will show in Section~\ref{ptio}, if 
condition (ii) in Theorem~\ref{ncsuf} is satisfied, then (i) holds for a 
subgroup $\,\Gm\subset\gp$. In addition, this $\,\Gm\hs$ may be chosen so that 
the compact quotient manifold $\,M=\hm\nnh/\hh\Gm$ is a bundle over the 
circle with some fibre $\,N\hs$ which is either a torus (when $\,\lz\,$ is a 
Lagrangian subspace of the solution space $\,\xe$, cf.\ Section~\ref{fols}), 
or a $\,2$-step nil\-man\-i\-fold (when $\,\lz\,$ is not Lagrangian).
\end{remark}
\begin{remark}\label{unimo}Since assertion (ii) in Theorem~\ref{ncsuf} 
includes the condition $\,\varPsi(\gkp)=\gkp$, it implies (\ref{pol}) for 
$\,P\,$ denoting the characteristic polynomials of both $\,\varPsi\,$ and 
$\,T:\lz\to\lz$. (Cf.\ the end of Section~\ref{prel}.) Thus, if (ii) holds, 
$\,T:\lz\to\lz\,$ must have determinant $\,\pm1$.
\end{remark}

\section{Nonexistence in dimension four}\label{nedf}
As mentioned in the Introduction, results of this section are used in the 
forthcoming paper \cite{derdzinski-roter}.

If $\,F:\mv\nh\to\mv$ is an endomorphism with the traceless part $\,E\,$ in 
a two\diml\ \rvs\ $\,\mv\nh$, then the traceless part of $\,F^2$ is 
$\,(\text{\rm tr}\,F)\hh E$. In fact, the matrix of $\,E\,$ in some basis is 
one of
\[
\left[\begin{matrix}\mu&\hskip2pt0\hskip-2pt\cr
0&-\hh\mu\end{matrix}\right],\hskip7pt
\left[\begin{matrix}0&-\hh\mu\cr
\mu&0\end{matrix}\right],\hskip7pt
\left[\begin{matrix}0&1\cr
0&0\end{matrix}\right],\hskip7pt
\]
with $\,\mu\in\bbR\hs$. Hence $\,E^2$ is a multiple of $\,1\,$ (the identity). 
The traceless and scalar parts of 
$\,F^2\nh=[E+(\text{\rm tr}\,F)/2\hs]^2\nh=
E^2\nh+(\text{\rm tr}\,F)\hh E\hs+\hs(\text{\rm tr}\,F)^2\nh/4$ thus are 
$\,(\text{\rm tr}\,F)\hh E\,$ and $\,E^2\nh+(\text{\rm tr}\,F)^2\nh/4$.
\begin{remark}\label{nozer}If $\,\dot\ro+\psi\ro=\delta\,$ for $\,C^1$ 
functions $\,\ro,\psi,\delta:\bbR\to\bbR\hs$, periodic with period $\,p>0$, 
and $\,\delta\ne0\,$ everywhere in $\,\bbR\hs$, then $\,\ro\ne0$ everywhere 
in $\,\bbR\hs$. In fact, the derivative $\,\dot\ro\,$ has the same nonzero 
signum at each zero of $\,\ro$, and so $\,\ro\,$ can have at most one zero 
in $\,\bbR\hs$, while $\,\ro$ with just one zero could not be periodic.
\end{remark}
\begin{lemma}\label{dfour}For $\,\fh,p,n,\mv\nh,\lr,A\,$ as in 
{\rm(\ref{obj})} with $\,n=4$, and a $\,C^\infty$ function 
$\,B:\bbR\to\hs\text{\rm End}\hs(\mv)$, periodic of period $\,p$, with 
$\,\dot B+B^2=\fh+A$, the determinant of the translation operator 
$\,T:\lz\to\lz$, defined at the beginning of Section~{\rm\ref{cecq}}, is not 
equal to $\,\pm1$.
\end{lemma}
\begin{proof}The traceless part of the equality $\,\dot B+B^2=\fh+A\,$ is 
$\,\dot E+(\text{\rm tr}\,B)\hh E=A$, where $\,E(t)\,$ denotes the 
traceless part of $\,B(t)$. (See the beginning of this section.) Choosing 
$\,v,v\hh'\in\mv\,$ with $\,\lg Av,v\hh'\rg\ne0\,$ and setting 
$\,\ro(t)=\lg E(t)v,v\hh'\rg$, $\,\hs\psi(t)=\hs\text{\rm tr}\,B(t)$, we now 
obtain $\,\dot\ro+\psi\ro=\delta$, where $\,\delta=\lg Av,v\hh'\rg\ne0$. 
Thus, $\,\ro\ne0\,$ everywhere by Remark~\ref{nozer}. Also, 
$\,\psi=(\delta-\dot\ro)/\nh\ro$, so that 
$\,\int_0^p\text{\rm tr}\,B(t)\,dt=\int_0^p\psi\,dt\ne0\,$ (as 
$\,\int_0^p(\dot\ro/\nh\ro)\,dt=0$). Our claim now follows from 
Lemma~\ref{rgsbs}(iii).
\end{proof}
\begin{theorem}\label{nonex}If\/ $\,\fh,p,\mv\nh,\lr,A\,$ are as in 
{\rm(\ref{obj})} for $\,n=4$, then, for $\,\hm,\gp\,$ defined in 
Section~\ref{csrr}, no subgroup $\,\Gm\hs$ of\/ $\,\gp\,$ acts on $\,\hm\,$ 
properly dis\-con\-tin\-u\-ous\-ly so as to produce a compact quotient 
manifold\/ $\,M=\hm\nnh/\hh\Gm$.
\end{theorem}
In fact, for such $\,\Gm\hs$ we might assume (\ref{smp}.a). Since (i) implies 
(ii) in Theorem~\ref{ncsuf}, it would follow from Remark~\ref{unimo} that the 
translation operator $\,T:\lz\to\lz\,$ has determinant $\,\pm1$, contrary to 
Lemma~\ref{dfour}.

\section{Proof that {\rm(ii)} implies {\rm(i)} in 
Theorem~\ref{ncsuf}}\label{ptio}
Given $\,B,\gkp,\varphi\,$ and $\,\rc\,$ with the properties listed in (ii), 
along with $\,\varPsi,\hh\varOmega\,$ as in (ii-b) and (ii-c), 
$\,\bbR\times\lz\,$ with the operation 
$\,(q,u)\cdot(r,w)=(q+r-\varOmega(u,w),\hs u+w)\,$ forms a Lie group, which 
we denote by $\,\hp$. In fact, the embedding 
$\,(q,u)\mapsto(0,q,u)\,$ identifies $\,\hp\,$ with a Lie subgroup of $\,\gp$, 
namely, $\,\Dt^{-1}(\lz)\,$ for $\,\Dt\,$ defined in (\ref{pgz}).

In view of (\ref{spc}.v) and (\ref{spc}.ii), $\,\hp\,$ is $\,2$-step 
nilpotent or Abelian, while $\,\varPsi:\hp\to\hp\,$ is a group automorphism 
(cf.\ (\ref{tao})). By (ii-c) and (ii-a), $\,\gkp\,$ is a subgroup of 
$\,\hp$, which gives rise to the quotient manifold 
$\,N=\hs\hp/\hh\gkp$, with $\,\gkp\,$ acting on $\,\hp\,$ by left 
translations.

Next, $\,N\hs$ is compact (a nil\-man\-i\-fold). In fact, suppose first that 
$\,\rc=0$. By (ii-c), $\,\varOmega(u,w)=0\,$ for all $\,u,w\in\lz$. (Namely, 
$\,\gkp\,$ spans $\,\bbR\times\lz$, and so $\,\Lm\,$ spans $\,\lz$.) Thus, 
$\,\hp\,$ is the additive group of $\,\bbR\times\lz$, and $\,N\hs$ is a torus. 
Now let $\,\rc>0$. Since $\,\gkp\,$ is a lattice in $\,\bbR\times\lz$, some 
compact set $\,Q\hs'\nh\subset\bbR\times\lz\,$ intersects every orbit of 
$\,\gkp\,$ (where $\,\gkp\,$ acts by translations in the additive group 
$\,\bbR\times\lz$), and so the image of $\,Q\hs'$ under the projection 
$\,\bbR\times\lz\to\lz\,$ is a compact set $\,\hq\subset\lz\,$ intersecting 
every orbit of the translation group $\,\Lm\subset\lz$. Compactness of 
$\,N\hs$ follows: the compact set 
$\,Q=[\hs0,\rc\hs]\times\hq\subset\hp=\bbR\times\lz\,$ intersects every orbit 
of $\,\gkp\,$ (now acting by left translations in $\,\hp$). To see this, 
consider the orbit of $\,(r,w)\in\hp$. The property of $\,\hq$, mentioned 
above, allows us to find $\,(q,u)\in\gkp\,$ with $\,u+w\in\hq$, and, for 
$\,l\in\bbZ\,$ such that $\,l\rc+q+r-\varOmega(u,w)\in[\hs0,\rc\hs]$, we get 
$\,(\rc,0)^l\nh\cdot(q,u)\cdot(r,w)\in Q\,$ (cf.\ (\ref{spc}.iii)), while 
$\,(\rc,0)\in\gkp\,$ by (ii-a).

Since $\,\varOmega:\xe\times\xe\to\bbR\,$ is nondegenerate (see the lines 
preceding (\ref{tao})), for some $\,u\in\xe\,$ the functional 
$\,-\hs2\hh\varOmega(u,\,\cdot\,)\,$ restricted to $\,\lz\,$ coincides with 
$\,\varphi$. Let us now fix any such $\,u\,$ and any 
$\,\widetilde q\in\bbR\hs$. For the unique $\,\widetilde u\in\lz\,$ with 
$\,\widetilde u(0)=u(0)$, (\ref{tao}) gives 
$\,\varOmega(T\hs\widetilde u,Tw)=\varOmega(\widetilde u,w)$, and so
\begin{equation}\label{map}
\hp\hs\in(r,w)\,\mapsto\,
(r+\widetilde q+\varOmega(\widetilde u-2u,w),T(w+\widetilde u\hs))\in\hs\hp
\end{equation}
describes the composite mapping in which the Lie-group automorphism 
$\,\varPsi:\hp\to\hp\,$ is followed by the right translation by the element 
$\,(\widetilde q,T\hs\widetilde u\hs)$ of $\,\hp=\bbR\times\lz$. As 
$\,\varPsi(\gkp)=\gkp\,$ (see (ii-b)) and $\,N=\hs\hp/\hh\gkp$, where 
$\,\gkp$ acts on $\,\hp\,$ by {\it left\/} translations, (\ref{map}) 
descends to the quotient nil\-man\-i\-fold $\,N=\hs\hp/\hh\gkp$, producing a 
\feo\ $\,\varPhi:N\nh\to N\nh$.

We now define $\,M\,$ to be the total space of a $\,C^\infty$ bundle over the 
circle $\,\bbR/\bbZ$, choosing our $\,(n-1)$\diml\ nil\-man\-i\-fold $\,N\hs$ 
to be the fibre, and using the \feo\ $\,\varPhi:N\nh\to N\hs$ to glue 
together the boundary components $\,\{0\}\times N\hs$ and $\,\{p\}\times N\hs$ 
of $\,[\hs0,p\hs]\times N\nh$, which we treat as copies of $\,N\nh$, while 
$\,\bbR/\bbZ\,$ is viewed as the result of identifying the two endpoints in 
$\,[\hs0,p\hs]$. Thus, $\,M=(\bbR\times N)/\bbZ$, where the action of 
$\,\bbZ\,$ on $\,\bbR\times N\hs$ is given by $\,k(t,y)=(t+kp,\varPhi^k(y))\,$ 
for $\,k\in\bbZ\,$ and $\,(t,y)\in\bbR\times N\nh$, or, equivalently, 
generated by the \feo\ $\,(t,y)\mapsto(t+p,\varPhi(y))$.

From now on we use the embedding $\,(q,u)\mapsto(0,q,u)\,$ to treat $\,\hp\,$ 
and $\,\gkp\,$ as subgroups of $\,\gp$. Let $\,\Gm\hs$ be the subgroup of 
$\,\hs\gp\,$ generated by $\,\gkp\,$ and the element $\,(1,q,u)$, with 
$\,q=\widetilde q\hs+\lg u(0),\dot u(0)-\nh B(0)\hh u(0)\rg\,$ and $\,u\,$ 
chosen above. Also, let $\,\pi:\hm\to M\,$ be the locally \feic\ surjective 
mapping arising as the composite 
$\,\hm\to\rto\nh\times\lz=\bbR\times\hp\to\bbR\times N\nh\to M$, in which the 
first arrow is the inverse of the equivariant \feo\ described in 
Remark~\ref{equiv}, the second sends $\,(t,(r,w))\,$ to $\,(t,\xi)$, where 
$\,\xi$ stands for the coset $\,\gkp\cdot(r,w)\,$ in $\,N=\hs\hp/\hh\gkp$, 
and the third one is the quotient projection 
$\,\bbR\times N\to M=(\bbR\times N)/\bbZ$.

It suffices to show that $\,\Gm\hs$ acts on $\,\hm\,$ freely and 
the $\,\pi$-pre\-im\-ages of points of $\,M\,$ coincide with the orbits of 
$\,\Gm$. (See Section~\ref{prel}.) To this end, first note that, by 
(\ref{spc}.iv) and (ii-b), $\,\gkp\,$ is invariant under the inner 
automorphism corresponding to $\,(1,q,u)$, and so any $\,(l,r,w)\in\Gm$, being 
a finite product of factors from the set $\,\gkp\cup\{(1,q,u),(1,q,u)^{-1}\}$, 
equals $\,(1,q,u)^k\nh\cdot(0,r\hs'\nnh,w\hh'\hh)\,$ for some 
$\,k\in\bbZ\,$ and $\,(0,r\hs'\nnh,w\hh'\hh)\in\gkp$. Secondly,
\begin{enumerate}
  \def\theenumi{{\rm\roman{enumi}}}
\item[($*$)] under the equivariant \feo\ $\,\hm\to\bbR\times\hp\,$ forming 
the first arrow in the composite $\,\pi:\hm\to M$, the actions (\ref{qku}.b) 
of $\,\hp\,$ and $\,\gkp\,$ on $\,\hm\,$ correspond to their actions on 
$\,\bbR\times\hp\,$ via left translations of the $\,\hp\,$ factor. (Cf.\ 
Remark~\ref{equiv} and (\ref{qkx}) with $\,k=0$.)
\end{enumerate}

The action of $\,\Gm\hs$ on $\,\hm\,$ is free: if $\,(l,r,w)\in\Gm\,$ has a 
fixed point, writing $\,(l,r,w)=(1,q,u)^k\nh\cdot(0,r\hs'\nnh,w\hh'\hh)\,$ as 
above, we obtain $\,\varPi(l,r,w)=k$, in view of (\ref{hom}); thus, $\,k=0\,$ 
by (\ref{qku}.b), and $\,(l,r,w)\in\gkp$, which shows that $\,(l,r,w)\,$ must 
be the identity, since $\,\gkp\,$ acts on $\,\hm\,$ freely (see ($*$)).

Also, two elements of the same orbit of $\,\Gm\hs$ in $\,\hm\,$ have the same 
image under $\,\pi$. In fact, \feo s $\,F:\hm\to\hm\,$ with 
$\,\pi\circ F=\pi\,$ form a group, which must contain $\,\Gm$, as it contains 
both $\,\gkp\,$ and $\,(1,q,u)$. (The former, by ($*$); the latter, since the 
action (\ref{qku}.b) of $\,(1,q,u)\,$ sending 
$\,\hmz=\{0\}\times\bbR\times\mv\hs$ onto $\,\hmp=\{p\}\times\bbR\times\mv\hs$ 
is easily verified to coincide with (\ref{map}) if one identifies $\,\hmz$ 
with $\,\{0\}\nh\times\hp\hs\approx\hs\hp\,$ and $\,\hmp$ with 
$\,\{p\}\times\hp\hs\approx\hs\hp\,$ via the \feo\ in ($*$), and uses the 
relation between $\,\lz\,$ and $\,B\,$ along with the definitions of $\,T\,$ 
and $\,\varOmega\,$ in Section~\ref{csrr}.)

Finally, let $\,(t,s,v),(t'\nnh,s\hh'\nnh,v\hs'\hh)\in\hm\,$ and 
$\,\pi(t,s,v)=\pi(t'\nnh,s\hh'\nnh,v\hs'\hh)$. Thus, $\,t'\nh=t-kp\,$ for some 
$\,k\in\bbZ$, and, replacing $\,(t'\nnh,s\hh'\nnh,v\hs'\hh)\,$ by the product 
$\,(1,q,u)^k\nh\cdot(t'\nnh,s\hh'\nnh,v\hs'\hh)$, we may assume that 
$\,t'\nh=t\,$ (cf.\ (\ref{qku}.b)). Then, by ($*$), $\,(t,s,v)\,$ and 
$\,(t'\nnh,s\hh'\nnh,v\hs'\hh)\,$ lie in the same orbit of $\,\gkp$.

\section{The main lemma}\label{male}
Given $\,p\in(0,\infty)$, let $\,\fvp$ be the set of all septuples 
$\,(\ax,\bx,\cx,\fh,\ay,\by,\cy)\,$ consisting of $\,C^\infty$ functions 
$\,\ax,\bx,\cx,\fh:\bbR\to\bbR\,$ of the variable $\,t$, periodic of 
period $\,p$, and constants $\,\ay,\by,\cy\in\bbR\,$ with $\,\ay+\by+\cy=0\,$ such 
that either $\,\by<\ay<\cy\,$ or $\,\by<\cy<\ay$, which satisfy, everywhere in 
$\,\bbR\hs$, the inequalities $\,\ax>\bx>\cx\,$ and the ordinary differential 
equations
\begin{equation}\label{thr}
\dot\ax+\ax^2=\fh+\hs\ay\hs,\hskip16pt\dot\bx+\bx^2=\fh+\hs\by\hs,\hskip16pt
\dot\cx+\cx^2=\fh+\hs\cy\hs.
\end{equation}
We denote by $\,\cv\,$ the subset of $\,\fvp$ formed by those 
$\,(\ax,\bx,\cx,\fh,\ay,\by,\cy)\,$ in which the functions 
$\,\ax,\bx,\cx\,$ and $\,\fh\,$ are all constant.
\begin{remark}\label{cnsta}If $\,(\ax,\bx,\cx,\fh,\ay,\by,\cy)\in\fvp$ and one 
of $\,\ax,\bx,\cx,\fh\,$ is constant, so are the other three. (In fact, a 
function $\,\ax\,$ such that $\,\dot\ax+\ax^2$ is constant cannot be periodic 
unless it is constant.)
\end{remark}
We now define a mapping 
$\,\hs\text{\rm spec}:\fvp\to\rtr$ by 
$\,\hs\text{\rm spec}\hs(\ax,\bx,\cx,\fh,\ay,\by,\cy)=(\lambda,\mu,\nu)\,$ 
for the unique $\,\lambda,\mu,\nu>0\,$ such that
\begin{equation}\label{spe}
\log\lambda\hs=\nh-\textstyle{\int_0^p\ax(t)\,dt}\hs,\hskip9pt
\log\mu\hs=\nh-\textstyle{\int_0^p\bx(t)\,dt}\hs,\hskip9pt
\log\nu\hs=\nh-\textstyle{\int_0^p\cx(t)\,dt}\hs.
\end{equation}
\begin{lemma}\label{surjc}The image 
$\,\hs\text{\rm spec}\hs(\fvp\nnh\smallsetminus\cv)\subset\rtr$ is the set of 
all\/ $\,(\lambda,\mu,\nu)\,$ satisfying conditions\/ {\rm(\ref{lmn})}.
\end{lemma}
We precede the proof of Lemma~\ref{surjc} with three other lemmas.
\begin{remark}\label{infdm}The most important part of Lemma~\ref{surjc} is the 
surjectivity claim, derived from a much stronger assertion. Namely, to conclude 
that the pre\-im\-age of every triple $\,(\lambda,\mu,\nu)\,$ with (\ref{lmn}) 
under the mapping  $\,\hs\text{\rm spec}:\fvp\nnh\smallsetminus\cv\to\rtr$ is 
nonempty, we show that it is, in fact, in\-fi\-nite\diml. See the final 
paragraph of this section.
\end{remark}
Any two out of the three equations (\ref{thr}) can be solved as follows.
\begin{lemma}\label{slveq}Let\/ $\,p\in(0,\infty)$, $\,\ay,\by\in\bbR\,$ and\/ 
$\,\ay\ne \by\hh$. Triples $\,(\ax,\bx,\fh)$ formed by $\,C^\infty$ 
functions $\,\bbR\to\bbR\hs$, periodic of period\/ $\,p$, with 
$\,\dot\ax+\ax^2=\fh+\ay\,$ and\/ $\,\dot\bx+\bx^2=\fh+\by$, are in a natural 
bijective correspondence with $\,C^\infty$ functions $\,\ro:\bbR\to\bbR\,$ 
which are periodic with period\/ $\,p\,$ and nonzero everywhere in 
$\,\bbR\hs$. The correspondence is given by $\,\ro=\ax-\bx$, and\/ $\,\ro\,$ 
determines $\,(\ax,\bx,\fh)\,$ via the relations 
$\,2\ax=\ro+(\ay-\by-\dot\ro)/\nh\ro$, 
$\,\hs2\bx=-\hs\ro+(\ay-\by-\dot\ro)/\nh\ro\,$ and\/ 
$\,\fh=\dot\ax+\ax^2-\ay$.
\end{lemma}
(This is a trivial exercise: for $\,\ro=\ax-\bx\,$ and $\,\psi=\ax+\bx$, the 
equations $\,\dot\ax+\ax^2=\fh+\ay$, $\,\,\dot\bx+\bx^2=\fh+\by\,$ give 
$\,\dot\ro+\psi\ro=\ay-\by$, so that $\,\ro\ne0$ everywhere by 
Remark~\ref{nozer}.) We will use the notation
\begin{equation}\label{drs}
\dy\,\,
=\,\,\{(\er,\es)\in\rto:\er>0\hh,\hskip6pt\es>0\hh,\hskip6pt\er\ne\es\}\hs.
\end{equation}
\begin{lemma}\label{sptoq}The mapping 
$\,(\ax,\bx,\cx,\fh,\ay,\by,\cy)\mapsto(\ro,\sy,\er,\es)$, with 
$\,\ro=\ax-\bx$, $\,\sy=\bx-\cx$, $\,\er=\ay-\by$, and\/ $\,\es=\cy-\by$, 
sends the set $\,\fvp$ defined in Section~{\rm\ref{male}} for any fixed\/ 
$\,p\in(0,\infty)\hs$ bijectively onto the set $\,\svp$ of all quadruples 
$\,(\ro,\sy,\er,\es)\,$ in which\/ $\,(\er,\es)\in\dy\,$ and\/ 
$\,\ro,\sy\,$ are positive $\,C^\infty$ functions of\/ $\,t\in\bbR\hs$, 
periodic with period\/ $\,p$, such that\/ 
$\,d\,[\hs\log\hs(\sy/\nh\ro)\hh]/dt\,=\,\ro\hs+\hs\sy\hs
-\hs\er\ro^{-1}\nh-\hs\es\hh\sy^{-1}\nnh$. The inverse mapping is given by 
$\,\ay=(2\er-\es)/3$, $\,\hs\by=-\hs(\er+\es)/3$, $\,\hs\cy=(2\es-\er)/3$, 
$\,\hs2\ax=\ro+(\er-\dot\ro)/\nh\ro$, 
$\,\hs2\bx=-\hs\ro+(\er-\dot\ro)/\nh\ro$, 
$\,\hs2\cx=-\hs\sy-(\es+\dot\sy)/\nh\sy\,$ and\/ $\,\fh=\dot\ax+\ax^2-\ay$.
\end{lemma}
In fact, the condition 
$\,d\,[\hs\log\hs(\sy/\nh\ro)\hh]/dt\,=\,\ro\hs+\hs\sy\hs-\hs\er\ro^{-1}\nh
-\hs\es\hh\sy^{-1}$ amounts to the equality between two expressions for 
$\,2\bx$, obtained by applying Lemma~\ref{slveq} first to $\,(\ax,\bx,\fh)\,$ 
and $\,\ay,\by$, then to $\,(\bx,\cx,\fh)\,$ and $\,\by,\cy$.
\begin{lemma}\label{sgphx}For any $\,p\in(0,\infty)$, the formula 
$\,x=\log\hs(\sy/\nh\ro)\,$ defines a bijective correspondence 
$\,(\ro,\sy,\er,\es)\mapsto(x,\er,\es)\,$ between the set $\,\svp$ defined in 
Lemma~{\rm\ref{sptoq}} and\/ $\,\xvp\nh\times\dy$, where $\,\xvp$ is the space 
of all\/ $\,C^\infty$ function $\,x\,$ of the variable $\,t\in\bbR\hs$, 
periodic with period\/ $\,p$, and\/ $\,\dy\,$ is the set\/ {\rm(\ref{drs})}. 
In terms of\/ $\,(x,\er,\es)$, the quadruple $\,(\ro,\sy,\er,\es)\in\svp$ is 
given by $\,2\hh\ro=(1+e^x)^{-1}
[\hs\dot x\hs+\hs\sqrt{\dot x^2+4(1+e^x)(\er+\es e^{-x})\,}\hs]\,$ and\/ 
$\,\sy=e^x\ro$. 
\end{lemma}
\begin{proof}If $\,d\,[\hs\log\hs(\sy/\nh\ro)\hh]/dt\,=\,\ro\hs+\hs\sy\hs
-\hs\er\ro^{-1}\nh-\hs\es\hh\sy^{-1}$ and $\,x=\log\hs(\sy/\nh\ro)$, replacing 
$\,\sy\,$ in the equality $\,\dot x=\ro+\sy-\er\ro^{-1}\nh-\es\hh\sy^{-1}$ by 
$\,e^x\nh\ro$, we get the required formula for $\,\ro\,$ by solving a 
quadratic equation.
\end{proof}
We now prove Lemma~\ref{surjc}. First, Lemmas~\ref{sptoq} and~\ref{sgphx} 
yield a bijection $\,\fvp\nh\to\xvp\nh\times\dy$, which we use to identify 
$\,\fvp$ and $\,\xvp\nh\times\dy$. This identification turns 
$\,\hs\text{\rm spec}\hs\,$ into a mapping 
$\,\xvp\nh\times\dy\to\rtr$ given by 
$\,\hs\text{\rm spec}\hs(x,\er,\es)=(\lambda,\mu,\nu)\,$ with 
$\,\lambda,\mu,\nu>0\,$ characterized by 
$\,(4\hh\log\lambda,\,4\hh\log\mu,\,4\hh\log\nu)
=(-\hs\dv-\ve,\,\dv-\ve,\,\dv-\ve+2\zt)$, where
\begin{equation}\label{dis}
\begin{array}{l}
\dv=\int_0^p(1+e^x)^{-1}
[\hs\dot x^2+4(1+e^x)(\er+\es e^{-x})\hs]^{1/2}\,dt\hs,\\
\ve=\int_0^p(1+\es e^{-x}\nnh/\nh\er)^{-1}
[\hs\dot x^2+4(1+e^x)(\er+\es e^{-x})\hs]^{1/2}\,dt\hs,\\
\zt=\int_0^p(1+e^{-x})^{-1}
[\hs\dot x^2+4(1+e^x)(\er+\es e^{-x})\hs]^{1/2}\,dt\hs.
\end{array}
\end{equation}
Thus, $\,\dv,\ve,\zt\,$ are all positive, $\,\ve\ne\zt\,$ (as $\,\es\ne\er$), and 
$\,\ve<\dv+\zt\,$ (since 
$\,(1+\es e^{-x}\nnh/\nh\er)^{-1}\nh<1=(1+e^x)^{-1}\nh+(1+e^{-x})^{-1}$). 
Consequently, $\,\lambda,\mu,\nu$ satisfy (\ref{lmn}). In other words, 
the values of $\,\hs\text{\rm spec}:\fvp\to\rtr$ all lie in 
the open set $\,\uv\,$ of all $\,(\lambda,\mu,\nu)\in\rtr$ with 
(\ref{lmn}).

Furthermore, $\,\hs\text{\rm spec}\hs\,$ restricted to 
$\,\cv\subset\fvp$ is a \feo\ $\,\cv\nh\to\hs\uv$. (This makes sense as 
$\,\cv\,$ is a submanifold of $\,\bbR^7\nnh$, namely, the graph of a 
$\,C^\infty$ mapping 
$\,\cv\hh'\ni(\ax,\bx,\cx)\mapsto(\fh,\ay,\by,\cy)\in\rfo\nnh$, where 
$\,\cv\hh'$ is the open set in $\,\rtr$ formed by all $\,(\ax,\bx,\cx)\,$ such 
that $\,|\bx\hh|<\ax\,$ and $\,|\bx\hh|<\nh-\hs\cx\ne\ax$.) In fact, by 
(\ref{spe}), $\,\hs\text{\rm spec}:\cv\to\uv\,$ is the composite 
$\,\cv\to\cv\hh'\nh\to\uv\,$ given by
\begin{equation}\label{cmp}
(\ax,\bx,\cx,\fh,\ay,\by,\cy)\mapsto(\ax,\bx,\cx)\mapsto(\lambda,\mu,\nu)
=(e^{-p\ax}\nnh,e^{-p\bx}\nnh,e^{-p\cx})\hs,
\end{equation}
and $\,(\ax,\bx,\cx)\mapsto(\lambda,\mu,\nu)\,$ in (\ref{cmp}) is a \feo\ 
$\,\cv\hh'\nh\to\uv$.

Finally, let $\,\yv\,$ be any vector subspace of $\,\xvp$ such that 
$\,2\le\dim\yv<\infty$ and $\,\yv\,$ contains the space of all constant 
functions (which we denote by $\,\bbR$). The subset $\,\yv\times\dy\,$ of 
$\,\xvp\nh\times\dy$ is at the same time an open set in $\,\yv\times\rto$, 
containing $\,\bbR\times\dy$. At every $\,(x,\er,\es)\in\bbR\times\dy$, 
the restriction of $\,\hs\text{\rm spec}:\xvp\nh\times\dy\to\uv\,$ to 
$\,\yv\times\dy\,$ is a submersion, since, as we just verified, its own 
restriction to $\,\bbR\times\dy\,$ is a \feo. (Under our identification 
$\,\fvp\approx\hs\xvp\nh\times\dy$, the subset $\,\cv\,$ of $\,\fvp$ 
corresponds to $\,\bbR\times\dy$, since, by Lemma~\ref{sgphx}, $\,x\,$ is 
constant if and only if both $\,\ro\,$ and $\,\sy\,$ are, which, in view of 
Lemma~\ref{sptoq}, amounts to constancy of $\,\ax,\bx,\cx\,$ and $\,\fh$.)

The pre\-im\-age of any given point $\,(\lambda,\mu,\nu)\in\uv\,$ under the 
mapping $\,\hs\text{\rm spec}:\xvp\nh\times\dy\to\uv\,$ thus contains a 
submanifold of dimension $\,\dim\yv-1$ intersecting $\,\bbR\times\dy\,$ at 
just one point, which completes the proof of Lemma~\ref{surjc}. The assertion 
about infinite dimensionality in Remark~\ref{infdm} now follows, since 
$\,\dim\yv\,$ can be arbitrarily large.

\section{Proof of Theorem~\ref{maith}}\label{pfmt}
We fix $\,p\in(0,\infty)\,$ as well as $\,k,l\in\bbZ\,$ with 
$\,2\le k<l\le k^2\nnh/4$, and choose the corresponding 
$\,\lambda,\mu,\nu\,$ with (\ref{lmn}) as in Lemma~\ref{polyn}. According to 
Lemma~\ref{surjc}, there exist $\,C^\infty$ functions 
$\,\ax,\bx,\cx,\fh:\bbR\to\bbR\hs$, periodic of period $\,p$, and constants 
$\,\ay,\by,\cy\in\bbR\,$ with $\,\ay+\by+\cy=0$, which satisfy equations 
(\ref{thr}) and (\ref{spe}), while $\,\fh\,$ is nonconstant 
(Remark~\ref{cnsta}), $\,\ax>\bx>\cx$, and $\,(\ay,\by,\cy)\ne(0,0,0)\,$ 
(since $\,\ay\ne\by\ne\cy\ne\ay$).

Let $\,\ao\in\hs\text{\rm End}\hs(\vo)\,$ and a $\,C^\infty$ function 
$\,\bo:\bbR\to\hs\text{\rm End}\hs(\vo)$, for a fixed $\,3$\diml\ \rvs\ 
$\,\vo$, be defined by requiring $\,\ao$ and all $\,\bo(t)\,$ to be 
simultaneously diagonalized by some fixed basis of $\,\vo$, with the 
eigenvalues $\,\ay,\by,\cy\,$ and, respectively, $\,\ax(t),\bx(t),\cx(t)$. 
By declaring the above basis orthonormal, we now introduce in $\,\vo$ a 
pseu\-\hbox{do\hs-}Euclid\-e\-an inner product of {\it arbitrary\/} 
signature.

Our $\,\fh,p,\vo,\ao$ and the inner product thus are objects of type 
(\ref{obj}) with $\,n=5$, while 
$\,\dot\bo\nh+\bo^2\nh=\fh+\ao$. By Lemma~\ref{rgsbs}, the space $\,\lo$ of 
all solutions $\,u:\bbR\to\vo$ to the equation $\,\dot u=\bo u\,$ is a 
Lagrangian subspace of the solution space $\,\eo$ for the equation 
$\,\ddot u=\fh u+\ao u$, and the the translation operator $\,T_0:\eo\to\eo$, 
given by $\,(T_0u)(t)=u(t-p)$, leaves $\,\lo$ invariant. According to 
Remark~\ref{cmmut}, $\,T_0:\lo\to\lo$ is diagonalizable with the eigenvalues 
$\,\lambda,\mu,\nu\,$ characterized by (\ref{spe}), so that $\,P\,$ appearing 
in Lemma~\ref{polyn} is its characteristic polynomial. As $\,P\,$ satisfies 
(\ref{pol}), there exists a lattice $\,\Lm_0$ in $\,\lo$ with 
$\,T_0(\Lm_0)=\Lm_0$. (See the end of Section~\ref{prel}.)

Next, we generalize this construction from $\,n=5\,$ to $\,n=3j+2$, for any 
integer $\,j\ge1$, using the original $\,\fh\,$ and $\,p$, but replacing each 
of $\,\vo,\ao,\bo,\eo,\lo,T_0$ and $\,\Lm_0$ by its $\,j\hh$th Cartesian power 
$\,\mv\nh,A,B,\xe,\lz,T$ and $\,\Lm$. Now $\,\mv\nh,A\,$ and each 
$\,B(t)\,$ is the direct sum of $\,j\,$ copies of $\,\vo,\ao$, or $\,\bo(t)$, 
so that $\,\xe,\lz\,$ and $\,T\,$ arise in the same way from $\,\eo,\lo$ and 
$\,T_0$. (We represent direct sums by Cartesian products.) Thus, 
$\,T(\Lm)=\Lm\,$ for the lattice $\,\Lm=\Lm_0\nh\times\ldots\times\Lm_0$ in 
$\,\lz=\lo\nh\times\ldots\times\lo$. In $\,\mv\hs$ we choose a 
pseu\-\hbox{do\hs-}Euclid\-e\-an inner product $\,\lr\,$ which is the 
orthogonal direct sum of inner products selected as above in each $\,\vo$ 
summand; the signature may vary from one summand to another, and so  
the resulting signature of $\,\lr\,$ is completely arbitrary.

The objects $\,\fh,p,\mv\nh,A\,$ and $\,\lr\,$ are again of type 
(\ref{obj}), this time with $\,n=3j+2$, and $\,\gkp=\bbZ\rc\times\Lm$, for 
any fixed $\,\rc\in(0,\infty)$, is a lattice in the vector space 
$\,\wv=\bbR\times\lz$. Our $\,B,\gkp\,$ and $\,\rc$, along with the zero 
functional $\,\varphi$, obviously satisfy condition (ii) in 
Theorem~\ref{ncsuf}. (Each $\hs B(t)\hs$ is diagonalized by an orthonormal 
basis of $\,\mv\nh$, so that $\,\lz\,$ is Lagrangian, 
cf.\ Lemma~\ref{rgsbs}(i), and $\,\varOmega(u,w)=0\,$ in (ii-c).)

Theorem~\ref{maith} is now immediate. Specifically, by Theorem~\ref{ncsuf}, 
there exists a subgroup $\,\Gm\subset\gp\,$ acting on $\,\hm\,$ freely and 
properly dis\-con\-tin\-u\-ous\-ly, such that the quotient manifold 
$\,M=\hm\nnh/\hh\Gm\,$ is compact. Lemmas~\ref{ecsrr} and~\ref{isome} imply 
that $\,M\,$ carries a metric $\,g\,$ with the properties required in 
Theorem~\ref{maith}. (The signature of the metric $\,\hg\,$ in (\ref{met}) 
is the result of augmenting, by one plus and one minus, the sign pattern of 
$\,\vh$, that is, of $\,\lr$.) The final clause of Theorem~\ref{maith} 
is in turn a consequence of Remark~\ref{nilbd} and Theorem~\ref{conse}(f), as 
$\,\Gm=\pi_1M$.
\begin{remark}\label{modul}The freedom of choosing 
$\,(\ax,\bx,\cx,\fh,\ay,\by,\cy)\,$ is in\-fi\-nite\diml\ 
(Remark~\ref{infdm}), which gives rise to an in\-fi\-nite\diml\ space of 
lo\-cal-i\-so\-met\-ry types of the resulting metrics $\,g$. In fact, the 
Ric\-ci tensor of the metric $\,\hg\,$ in (\ref{met}) is a constant 
multiple of $\,\fh(t)\,dt\otimes dt$, and the $\,1$-form $\,dt\,$ is 
$\,\hg$\prl\ \cite[p.\ 93]{roter}. Therefore, the function $\,\fh$ constitutes 
a local geometric invariant of $\,\hg$, defined up to affine changes of the 
variable $\,t\,$ and multiplications of $\,\fh\,$ by nonzero constants.
\end{remark}

\end{document}